\documentclass[lettersize,journal]{IEEEtran}

\usepackage{amsmath,amsfonts}
\usepackage{algorithmic}
\usepackage{algorithm}
\usepackage{array}
\usepackage[caption=false,font=normalsize,labelfont=sf,textfont=sf]{subfig}
\usepackage{textcomp}
\usepackage{stfloats}
\usepackage{url}
\usepackage{verbatim}
\usepackage{graphicx}
\usepackage{xcolor}
\usepackage{cite}
\usepackage{multirow}
\usepackage{ifthen}
\newboolean{arxiv}

\newtheorem{theorem}{Theorem}

\newtheorem{definition}{Definition}
\newtheorem{corollary}{Corollary}

\newtheorem{assumption}{Rule}

\newcounter{boxlblcounter}  
\newenvironment{boxlabel}
  {\begin{list}
    {\arabic{boxlblcounter}}
    {\usecounter{boxlblcounter}
     \setlength{\labelsep}{0.25em}
     \setlength{\itemindent}{0em} 
     \setlength{\leftmargin}{0.15cm}
     \setlength{\rightmargin}{0 cm}
     
    }
  }
{\end{list}}

\begin{document}

\setboolean{arxiv}{true}  

\title{Market Power Mitigation in Two-stage Electricity Markets with Supply Function and Quantity Bidding}
% Two-stage Electricity Market Equilibrium under Generator Supply Function Bidding with Market Power Mitigation and Strategic Load

\author{Rajni Kant Bansal, Yue Chen, Pengcheng You, Enrique Mallada
        % <-this % stops a space
\thanks{R. K. Bansal and E. Mallada are with the Johns Hopkins University, Baltimore, MD 21218, USA {(email: \{rbansal3,mallada\}@jhu.edu)}}% <-this % stops a space
%\thanks{Manuscript received April 19, 2021; revised August 16, 2021.}

\thanks{Y. Chen is with The Chinese University of Hong Kong, Shatin, Hong Kong SAR, China {(email:yuechen@mae.cuhk.edu.hk)}}

\thanks{P. You is with the Peking University, Beijing, China  {(email:pcyou@pku.edu.cn)}}

\thanks{This work was supported by NSFC (grants  72201007, T2121002, 72131001, 61973163, and 52307144) and NSF (grants CAREER ECCS 1752362, CPS ECCS
2136324, and Global Centers 2330450).}
}



\maketitle

\begin{abstract}
%v8
\textcolor{black}{
%v1_revised
Two-stage settlement electricity markets, which include day-ahead and real-time markets, often observe undesirable price manipulation due to the price difference across stages, inadequate competition, and unforeseen circumstances. To mitigate this, some Independent System Operators (ISOs) have proposed system-level market power mitigation (MPM) policies in addition to existing local policies. These system-level policies aim to substitute noncompetitive bids with a default bid based on estimated generator costs. However, without accounting for the conflicting interest of participants, they may lead to unintended consequences when implemented. In this paper, we model the competition between generators (bidding supply functions) and loads (bidding quantity) in a two-stage market with a stage-wise MPM policy. An equilibrium analysis shows that a real-time MPM policy leads to equilibrium loss, meaning no stable market outcome (Nash equilibrium) exists. A day-ahead MPM policy leads to Stackelberg-Nash game, with loads acting as leaders and generators as followers. Despite estimation errors, the competitive equilibrium is efficient, while the Nash equilibrium is comparatively robust to price manipulations. Moreover, analysis of inelastic loads shows their tendency to shift allocation and manipulate prices in the market. Numerical studies illustrate the impact of cost estimation errors, heterogeneity in generation cost, and load size on market equilibrium.
}
\end{abstract}


\section{Introduction}
\IEEEPARstart{M}{ost} wholesale energy markets in the US consider a two-stage settlement system as a market norm, i.e., day-ahead and real-time markets. The first stage, the day-ahead (forward) market, clears a day before the delivery based on the hourly forecasts of resources for the next day and accounts for the majority of energy trades. The second stage, the real-time (spot) market, occurs at a faster timescale (typically every five minutes) and is considered a last resort for participants to adjust their commitment following forecast errors~\cite{pyou_discovery,bansal_e_energy}. %two_stage,two_stage_design}. 
The main goal of such a sequential two-stage market is to operate efficiently and encourage market participation. However, the often price difference between the two stages in practice, due to intrinsic uncertainty in the forecast, unscheduled maintenance, etc., creates opportunities for price speculation and arbitrage, which could be further exploited by strategic participants to their benefit~\cite{low_lmp_spread,borensteinpricespread2008}. %signaling efficiency losses~\cite{low_lmp_spread,borensteinpricespread2008}. %,market_power_strat_gen_adam}. 

To discourage suppliers from exploiting consumers, most operators employ an inbuilt local market power mitigation mechanism (LMPM) triggered at congestion during market clearing~\cite{local_mpm_stanford, caiso_day_ahead}. Despite this, some operators, like California Independent System Operator (CAISO), have documented periods of time with non-competitive bids (approximately $2\%$ hours in the case of CAISO~\cite{caiso_scop}). It led to the development of initiatives aimed at implementing \emph{system-level} market power mitigation (MPM), i.e., bid mitigation similar to LMPM, but system-wide for each stage separately~\cite{caiso_proposal,caiso_extda-bundle1}. Such system-level policies, when implemented, substitute in, e.g., real-time or day-ahead, any non-competitive bids with \emph{default bids}, which estimate generator costs based on the operator’s knowledge of technology, fuel prices, and operational constraints\color{black}~\cite[\S 39.7.1]{caiso_tarif},\cite{caiso_bid_memo}\color{black}. Although such market policies are straightforward, their effect on market outcome remains unknown if implemented without accounting for the conflicting interest of individual participants. This paper studies the proposed system-level policies and discusses the possible unintended effects.

%Even though it is reasonable but too aggressive to implement such mitigation policies in both day-ahead and real-time, these market modifications are being considered separately, starting with the real-time stage in the first phase and followed by the day-ahead. In particular, CAISO argues that real-time is more susceptible to market power and day-ahead is relatively competitive due to additional measures like virtual bidding that add competitive pressure on market clearing~\cite{caiso_proposal}.  

Precisely, we study a sequential game formulation in a two-stage market with an MPM policy to analyze the competition between generators (bidding supply functions and seeking to maximize individual profit)~\cite{Joh2004Efficiency} and loads (bidding demand quantities and minimizing payment)~\cite{rad_fixed_bid} such that the market operator substitutes generators' bids with default bids as per the policy. \textcolor{black}{In this paper, we assume that an operator makes an error in estimating the truthful cost of dispatching generators in a stage with an MPM policy.}
We show that a real-time MPM policy results in a loss of market equilibrium. However, the complimentary case of a day-ahead MPM policy leads to a form of \emph{Stackelberg-Nash game} with loads leading generators in their decision-making. \textcolor{black}{A detailed Nash equilibrium analysis for this case shows a stable market outcome that is comparatively robust to price manipulations.}

The main contributions of this paper are summarized below:

\begin{boxlabel}
    \item We show that a real-time MPM policy leads to a Nash game in the day-ahead, while generators participate truthfully in real-time. \textcolor{black}{We characterize the competitive equilibrium of such a game, which is inefficient w.r.t the social planner’s problem}. Further, competition between price-anticipating participants does not result in a stable market outcome, and a Nash equilibrium does not exist. 

    \item We then study the impact of a day-ahead MPM policy that leads to a generalized Stackelberg-Nash game with loads acting as leaders in the day-ahead market and generators acting as followers in the real-time market. \textcolor{black}{Despite the operator's error in cost estimation, the competitive equilibrium of the resulting game is efficient. Also, the Nash equilibrium, assuming that generators are homogeneous and bid symmetrically for closed-form analysis, is robust to price manipulations compared to standard markets, i.e., a two-stage market without any mitigation policies.} 

    \item To understand the impact of these policies, we compare the day-ahead MPM policy market equilibrium with the equilibrium in a standard market. The closed-form analysis shows that prices across stages are the same for the two cases. However, loads acting as leaders in a market with a day-ahead MPM policy allocates higher demand in day-ahead at the expense of generators' profit. \textcolor{black}{Despite being inelastic, loads can shift their allocation and manipulate prices in the market.}
    %\textcolor{red}{It leads to an unintended consequence of additional market power instead of mitigation compared to the standard markets.}
    %\enrique{In which sense the market power is additional?}
    
    %\item  To understand the strategic behavior of players at the Nash equilibrium in a market with a day-ahead MPM policy, we further characterize (artificial) unilateral Nash equilibrium, i.e., either loads or generators act strategically. Our analysis shows that the generators, even when participating strategically in the market, exercise limited market power and loads acting as a leader manages to influence the market in their favor. 

    \item We further provide a detailed numerical study to illustrate the impact of a day-ahead MPM policy. We show that the Nash equilibrium converges to competitive equilibrium as the number of participants increases in the market. \textcolor{black}{Overestimation of the cost of generators benefits them with a higher profit and helps mitigate the market power of loads.} Furthermore, the case with heterogeneity in generator cost shows that expensive generators are least affected when benchmarked with the competitive equilibrium. The case with significant diversity in load participants reveals that a sufficiently smaller load could earn money instead of making payments at the expense of larger loads in the market.

\end{boxlabel}

\textit{Related work:} 
% revised_v1
%To date, the problem of characterizing the effect of \emph{system-level} market power mitigation policy based on default bids in two-stage markets has not been considered by the literature.Numerous works have developed ways to identify existence of market power~\cite{survey_all_1,survey_all_2}, characterized equilibrium games to investigate competition between various market players~\cite{two_stage_competition_1,bansal2021market}, and analysed mechanisms for market power mitigation~\cite{yue_mpm,mpm_mechanism_1,mpm_mechanism_2,mpm_mechanism_3,cai_forward_contract}. In particular, the role of forward contracting~\cite{cai_forward_contract}, demand shifting~\cite{mpm_mechanism_3}, bidding capacity division and constraint~\cite{mpm_mechanism_2}, network expansion~\cite{mpm_mechanism_1}, promoting interconnection and capacity regulation~\cite{newbery2002mitigating} and LMP-based default energy bids~\cite{DEB_MPM}, as a tool for local market power mitigation.  Despite the extensive studies on the cross-group market competition with only a few papers conducting counterfactual analysis, i.e., investigating the effect of an MPM policy on the market, to the best of our knowledge, our work is the first one to study the impact of a system-level market power mitigation policy based on default bids in a two-stage market and characterize the behavior of inelastic demand on the market. Our work formally analyzes the effect of default bid substitution on the market outcome, provides theoretical guarantees, and discusses recommendations to policymakers.
%
% revised_v2
\textcolor{black}{
Understanding market power and strategies for mitigating it has been an extensive subject of study in the literature. Prior works have studied the identification of market power~\cite{survey_all_1,survey_all_2}, the development of metrics to quantify it~\cite{mkt_pow_new_id}, competition between various market players\cite{two_stage_competition_1,bansal2021market,yue_mpm}, and general analysis of market power in two-stage markets~\cite{pyou_discovery,mpm_mechanism_1,mpm_mechanism_2,mpm_mechanism_3,cai_forward_contract}. In addition, some works have analyzed local MPM policies, e.g. see~\cite{lmpm_counter_factual}. However, a closed-form analysis, where the system-level policy effect on the resulting equilibrium is studied counterfactually, is rare in the literature. Our work contributes to the field by conducting a counterfactual study on the impact of CAISO's system-level policy based on default bids and the role of inelastic demand in the market. To the best of our knowledge, there is no existing study that investigates these features.
%and counterfactual analysis of role of mitigation mechanisms~\cite{}. In particular the role of forward contracting~\cite{cai_forward_contract}, demand shifting~\cite{mpm_mechanism_3}, bidding capacity division and constraint~\cite{mpm_mechanism_2}, network expansion~\cite{mpm_mechanism_1}, interconnection and capacity regulation~\cite{newbery2002mitigating} and LMP-based default energy bids~\cite{DEB_MPM}, for local market power mitigation. These attempts to study market power have improved our understanding of the field and support ISOs' attempts to restrict market power with a system of checks and balances. However, these policies are not generalizable, and market power still exists due to the lack of a system-level initiative. Our work extends this literature by conducting a counterfactual study on the impact of system-level policy based on default bids. It provides a closed-form analysis of such a policy, characterizes the behavior of inelastic loads with generators in a two-stage settlement market, provides theoretical guarantees, and discusses recommendations to policymakers. To the best of our knowledge, there is no existing study that investigates these features.
}

\textit{Paper Organization:}
The rest of the paper is structured as follows. In Section~\ref{sec_2} we first introduce the social planner problem, two-stage market model, and participants' behavior, and then define a two-stage market equilibrium. In Section~\ref{sec_3} we model the market power mitigation policy for each stage, characterize the market equilibrium, and compare it with the solution to the social planner problem. We first compare the impact of MPM policies on market equilibrium and then compare it with a standard market equilibrium in Section~\ref{sec_4}. Numerical studies on market power for a day-ahead MPM policy, limitations of work along with policy implications, and conclusions are in Section~\ref{sec_5}, \ref{sec_pre_6}, and ~\ref{sec_6}, respectively.

\subsubsection*{Notation} 
\textcolor{black}{We use standard notation $f(a,b)$ to denote a function of independent variables $a$ and $b$. However, we use $f(a;b)$ to represent a function of independent variable $a$ and parameter $b$.}  

\section{Market Model}\label{sec_2}

In this section, we formulate the social planner problem and describe the standard two-stage settlement market. We then formally define participants' behavior, i.e., price-taking or price-anticipating, and lay out a general market equilibrium. 

\subsection{Social Planner Problem}

Consider a single-interval two-stage settlement market where a set $\mathcal{G}$ of generators participate with a set $\mathcal{L}$ of inelastic loads to meet inelastic aggregate demand $d \in \mathbb{R}$. Each generator $j \in \mathcal{G}$ supplies $g_j \in  \mathbb{R}$ and each inelastic load $l\in\mathcal{L}$ consumes $d_l \in \mathbb{R}$ respectively, where $\sum_{l\in\mathcal{L}}d_l = d$. We define $G:= |\mathcal{G}|$ and $L:= |\mathcal{L}|$ to denote the number of generators and loads, respectively. Assuming a quadratic cost function for the generators, the social planner problem  --- minimum cost of meeting aggregate inelastic demand --- is given by 
\begin{subequations}
\label{planner_problem}%
\begin{eqnarray}
    \min_{g_j,j\in\mathcal{G}} && \sum\nolimits_{j\in \mathcal{G}} \frac{c_j}{2}{g_j}^2 \label{planner_obj} \\
    \text{s.t.} \!\!\!\!
    && \sum\nolimits_{l\in\mathcal{L}}d_l =\sum\nolimits_{j\in\mathcal{G}}g_j \label{power_bal_constraint}
\end{eqnarray}
\end{subequations}
where \eqref{power_bal_constraint} enforces the power balance in the market. 

\subsection{Two-Stage Market Mechanism}

\begin{figure}
    \centering
    \includegraphics[width=0.75\linewidth]{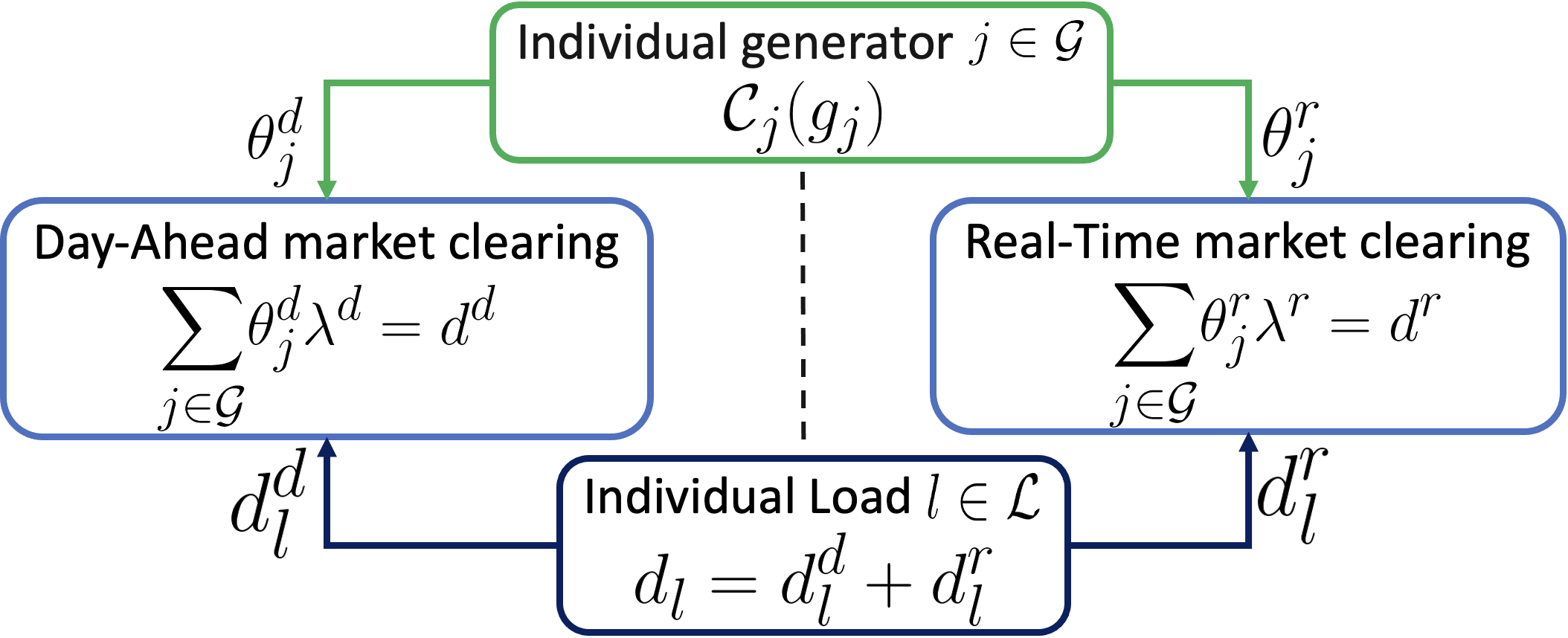}
    \caption{Two-stage Market Mechanism}%\emmargin{Please update the figure.}
    \label{fig:two_stage}
\end{figure}

In this subsection, we define the two-stage market clearing, as shown in Figure~\ref{fig:two_stage}. \textcolor{black}{The net output $g_j$ of each generator $j$ and individual load $d_l$ of load $l$ is allocated over two stages, such that
\begin{align}\label{real_time_load_condition}
    g_j = g_j^d+g_j^r, \ d_l^d +d_l^r = d_l
\end{align} 
where $(g_j^d, d_l^d)$ and $(g_j^r,d_l^r)$ represent allocation in day-ahead and real-time markets, respectively. %Similarly, each load $l$ allocates its inelastic demand $d_l$ over two stages %between $d_l^d$ and $d_l^r$ 
%such that
%\begin{align} \label{real_time_load_condition}
%    d_l^d +d_l^r = d_l
%\end{align}
%where $d_l^d, d_l^r$ represent allocation in DA and RT, respectively. 
}
\subsubsection{Day-Ahead Market} The power output of each generator $j \in \mathcal{G}$ in the day-ahead market is denoted by $g_j^{d}$. Each generator $j $ submits a supply function parameterized by the slope $\theta_j^d$, that indicates willingness of generator $j$ to supply $g_j^{d}$ as a function of price
\begin{align}
    g_j^{d} =\theta_j^d\lambda^{d} \label{gen_da_bid}
\end{align}
where $\lambda^{d}$ denotes the price in the day-ahead market. Each load $l\in\mathcal{L}$ bids quantity $d_l^{d}$ in the day-ahead market. Based on the bids $(\theta_j^{d}, d_l^{d})$ from participants, the market operator clears the day-ahead market to meet the supply-demand balance.
\begin{align}%\label{da_dispatch}
    %&\sum_{j\in\mathcal{G}} g_j^{d} = \sum_{l\in\mathcal{L}} d_l^{d} = d^{d}  \nonumber \\
     &  \sum\nolimits_{j\in\mathcal{G}} \theta_j^d\lambda^{d} = d^{d}  \label{da_power_bal}
\end{align}
The optimal solution to the day-ahead dispatch problem~\eqref{da_power_bal} gives the optimal dispatch $(g_j^{d},d_l^{d})$ and clearing prices $\lambda^{d}$ to all the participants. Each generator $j\in\mathcal{G}$ and load $l\in\mathcal{L}$ are paid $\lambda^{d}g_j^{d}$ and $\lambda^{d}d_l^{d}$ as part of the market settlement.

\subsubsection{Real-Time Market}

The power output of each generator $j$ in real-time market is denoted by $g_j^{r}$ and their bid is:
\begin{align}\label{gen_rt_bid}
    g_j^{r} = \theta_j^{r}\lambda^{r}
\end{align}
where $\lambda^{r}$ denotes the price in the real-time market. The supply function bid is parameterized by $\theta_j^{r}$, indicating willingness of generator $j$ to supply $g_j^{r}$ at the price $\lambda^{r}$. Each load $l\in\mathcal L$ submits quantity bids $d_l^{r}$. Given the bids $(\theta_j^{r}, d_l^{r})$, the operator clears the real-time market to meet the supply-demand balance.
\begin{align}%\label{rt_dispatch}
    %& \sum_{j\in\mathcal{G}} g_j^{r} = \sum_{l\in\mathcal{L}} d_l^{r} = d^{r} \nonumber \\
    & \sum\nolimits_{j\in\mathcal{G}} \theta_j^{r}\lambda^{r} = d^{r}  \label{rt_power_bal}
\end{align}
Similar to the day-ahead market clearing, the optimal solution to the dispatch problem~\eqref{rt_power_bal} gives the optimal dispatch and the market clearing prices $\lambda^{r}$ to all the participants, such that each generator $j\in\mathcal{G}$ and load $l\in\mathcal{L}$ produces or consumes $g_j^{r}$ and $d_l^{r}$, and is paid or charged $\lambda^{r}g_j^{r}$ and $\lambda^{r}d_l^{r}$, respectively. 

\subsubsection{Market Rules and Goal}
In this section, we first define a set of rules to account for degenerate cases in the market mechanism and then discuss the goal of a two-stage market.  
\begin{assumption}\label{assumpt_same_prc}
    For $v\in \{d,r\}$ and $w\in \{d,r\}$, if the net supply and demand of the generators and loads in a stage follow 
    \begin{align}
        \textcolor{black}{\sum\nolimits_{j\in\mathcal{G}} \theta_j^v\lambda^{v} = 0, \ d^v = 0 \implies \lambda^v = \lambda^w,  \ v \neq w} %\ v \subset \{d,r\}
    \end{align}
    i.e., the clearing price in that stage is set to the clearing prices of the other stage with a non-zero demand.
\end{assumption}
\begin{assumption}\label{assumpt_zero_prc}
    For $v\in \{d,r\}$, if the net supply and net demand of the generators and loads in a stage follow 
    \begin{align}
        \sum\nolimits_{j\in\mathcal{G}}\theta_j^v\lambda^{v} = 0, \ d^v \neq 0 \implies \lambda^v = 0 %\ v \subset \{d,r\}
    \end{align}
    i.e., the clearing price is set to zero, and demand is split evenly across all the loads.
\end{assumption}
We are interested in two-stage market outcomes that satisfy 
%\begin{subequations}
\begin{eqnarray}
    %& g_j^{d}+g_j^{r} = g_j, \ j \in \mathcal{G} \label{two_stage_gen_condition}\\
    %& d_l^{d}+d_l^{r} = d_l, \ l \in \mathcal{L} \label{real_time_load_condition} \\
    & \sum\limits_{j\in\mathcal{G}} (g_j^{d}+g_j^{r}) = \sum\limits_{j\in\mathcal{G}} g_j = \sum\limits_{l \in \mathcal{L}} (d_l^{d}+d_l^{r})= \sum\limits_{l \in \mathcal{L}}d_l = d \label{two_stage_power_bal}
\end{eqnarray}
%\end{subequations}
and solve the social planner problem~\eqref{planner_problem}. Though the market outcome may deviate from the optimal social planner solution, signaling efficiency losses due to price manipulation by participants, we quantify such deviations to understand the behavior of participants and the market outcome.

\subsection{Participant Behaviour}

In this section, for the purposes of our study, we introduce two different types of rational participants' behavior, price-taking, and price-anticipating. Each generator $j\in\mathcal{G}$ seeks to maximize their profit $\pi_j$, given by: 
%\noindent\emph{Generator Profit}
\begin{align}\label{generator_profit}
    &\pi_j(g_j^{d},g_j^{r},\lambda^{d},\lambda^{r}) \!:= \! \lambda^{r}g_j^{r} \!\!+ \lambda^{d}g_j^{d} \!-\! \frac{c_j}{2}{g_j}^2
\end{align}
%Substituting the supply function~\eqref{gen_da_bid} and \eqref{gen_rt_bid} in~\eqref{generator_profit}, we get 
%\begin{align}\label{generator_profit}
%    &\!\!\!\pi_j(\theta_j^{d},\!\theta_j^{r},\!\lambda^{d}\!,\!\lambda^{r}) \! =\!  \lambda^{d}h(\lambda^{d};\theta_j^{d}) \!+\! \lambda^{r}h(\lambda^{r};\theta_j^{r})\!-\! \mathcal{C}_j\!\left(\Tilde{h}(\theta_j^{d},\!\theta_j^{r},\!\lambda^{d}\!,\!\lambda^{r})\right)\label{generator_profit_bids}
%\end{align}
%where $\Tilde{h}$
Each load $l \in \mathcal{L}$ aims to minimize their payments $\rho_l$, as:
%\noindent\emph{Load Payment}
\begin{align}\label{load_payment_definition}
    & \rho_l(d_l^{d},d_l^{r},\lambda^{d},\lambda^{r}) :=  \lambda^{d}d_l^{d} + \lambda^{r}d_l^{r}
\end{align}
Substituting the load coupling constraint~\eqref{real_time_load_condition} in~\eqref{load_payment_definition} we get,
\begin{align}%\label{load_payment_intermediate_rt_mpm}
    & \rho_l(d_l^{d},\lambda^{d},\lambda^{r}) :=  \lambda^{d}d_l^{d} + \lambda^{r}(d_l-d_l^{d}) \label{load_payment_bids}
\end{align}
For each load $l \in \mathcal{L}$, the allocation in the day-ahead market $d_l^{d}$ determines its allocation in the real-time market $d_l^{r}$ due to the demand inelasticity.

\subsubsection{Price-Taking Participants}

A price-taker participant is defined below:
\begin{definition}
    A market participant is price-taking if it accepts the existing prices in the market and does not anticipate the impact of its bid on the market prices.
\end{definition}

Given the prices in the day-ahead market $\lambda^{d}$ and real-time market $\lambda^{r}$, the generator individual problem is given by:
%\noindent\emph{Price-taking Generator Bidding problem}
\begin{align}\label{generator_price_taking_profit}
    &\max_{g_j^{d},g_j^{r}} \  \pi_j(g_j^{d},g_j^{r};\lambda^{d},\lambda^{r}) 
\end{align}
Similarly given the prices $\lambda^{d}, \lambda^{r}$, the individual bidding problem for load is given by:
%\noindent\emph{Price-taking Load Bidding problem}
\begin{align}\label{load_price_taking_payment}
    \min_{d_l^{d}} \  &  \rho_l(d_l^{d};\lambda^{d},\lambda^{r})
\end{align}
 We next define the price-anticipating (or strategic) participants.

\subsubsection{Price-Anticipating Participants}

A price-anticipating participant is defined below:

\begin{definition}
    A market participant is price-anticipating (strategic) if it anticipates the impact of its bid on the prices in two stages and has complete knowledge of other participants' bids. 
\end{definition}

The individual problem of a price-anticipating generator is:
%\noindent\emph{Strategic Generator Bidding problem}
\begin{subequations}\label{generator_strategic_profit_total}
\begin{align}\label{generator_strategic_profit}
    & \!\!\!\max_{g_j^{d},g_j^{r},\lambda^d,
    \lambda^r} \ \! \pi_j \! \left(g_j^{d},g_j^{r},\lambda^{d}\!\left(g_j^{d};\overline{g}_{-j}^{d},d^{d}\right),\lambda^{r}\!\!\left(g_j^{r};\overline{g}_{-j}^{r},d^{r}\right) \right) \\
    &\textrm{ s.t. } \eqref{da_power_bal}, \eqref{rt_power_bal}
\end{align}
\end{subequations}
where $\overline{g}_{-j}^{d} := \sum_{k \in \mathcal{G}, k \neq j}g_k^{d}$, and $\overline{g}_{-j}^{r} := \sum_{k \in \mathcal{G}, k \neq j}g_k^{r}$ .
The generator $j$ maximizes its profit while anticipating the market clearing prices in the day-ahead and real-time market~\eqref{da_power_bal},\eqref{rt_power_bal}, along with complete knowledge of load bids $d_l^{d}, d_l^{r}, l\in\mathcal{L}$, and other generators' bids $\theta_k^d,\theta_k^r, k \in \mathcal{G}, k \neq j$. Similarly, the individual problem for strategic load $l$ with complete knowledge of prices in two stages~\eqref{da_power_bal},\eqref{rt_power_bal} and other participants' bids:
%\noindent\emph{Strategic Load Bidding problem}
\begin{subequations}\label{load_strategic_payment}
\begin{align}
    & \!\!\!\min_{d_l^{d},\lambda^d,\lambda^r} \ \rho_l\!\left(d_l^{d}, \lambda^{d}\!\left(d_l^d;g_j^{d},\overline{d}_{-l}^{d}\right)\!\!,\lambda^{r}\!\!\left(d_l^d;g_j^{r},\overline{d}_{-l}^{r}\right)\right)\\
    &\textrm{ s.t. } \eqref{da_power_bal}, \eqref{rt_power_bal}
\end{align}
\end{subequations}
where the load $l$ minimizes its payment in the market and $\overline{d}_{-l}^{d} := \sum_{l \in \mathcal{L}, k \neq l}d_l^{d}, \ \overline{d}_{-l}^{r} := \sum_{l \in \mathcal{L}, k \neq l}d_l^{r}$.

\subsection{Market Equilibrium}

In this section, for the purpose of this study, we characterize the market equilibrium in a two-stage settlement electricity market. At the equilibrium, no participant has any incentive to deviate from their bid, and market clears, as defined below.
%}
\begin{definition} \label{market_eqbm}
We say the participant bids and market clearing prices $( {\theta_j^{d}}, {\theta_j^{r}},j\in \mathcal{G}, d_l^{d}, d_l^{r}, l \in\mathcal{L},\lambda^{d}, \lambda^{r})$ in the day-ahead and real-time respectively form a two-stage market equilibrium if the following conditions are satisfied:
\begin{enumerate}
    \item For each generator $j \in \mathcal{G}$, the bid $\theta_j^{d}, \theta_j^{r}$ maximizes their individual profit.
    \item For each load $l \in \mathcal{L}$, the allocation $d_l^{d}, d_l^{r}$ minimizes their individual payment.
    \item The inelastic demand $d \in \mathbb{R}$ is satisfied with the market-clearing prices $\lambda^{d}$ given by~\eqref{da_power_bal} and $\lambda^{r}$ given by~\eqref{rt_power_bal} over the two-stages of the market.
\end{enumerate}
\end{definition}

We will study market equilibria as a tool to understand the impact of MPM policies.

\section{Understanding Impact of MPM Policy}\label{sec_3}

In this section, we first characterize the market equilibrium in a standard two-stage market without any such mitigation policy~\cite{pyou_discovery}, then model mitigation policies and the resulting market equilibrium. In particular, ISOs have significant prior knowledge of market participants allowing them to evaluate the competitiveness of energy bids. For example, operators are aware of the generator’s technology, fuel prices, and operational constraints that can be used to estimate or bound the generator’s cost~\cite[\S 39.7.1]{caiso_tarif},\cite{caiso_bid_memo} within a reasonable threshold under the mitigation policies. \textcolor{black}{We assume that the operator makes an error in estimating the truthful cost of dispatching the generator in a stage with a mitigation policy.} 
%Also, we assume a quadratic cost function for each generator $j$, parameterized by cost coefficient $c_j$, 
%\begin{align}
%    \mathcal{C}_j(g_j) = \frac{c_j}{2}{g_j}^2 \label{gen_quad_cost}
%\end{align}
%We assume that the supply function bid $h(\lambda^{v};\theta_j^{v})$ of each generator $j$ is given by:
%\begin{align}
%    g_j^{v} := h(\lambda^{v};\theta_j^{v}) = \theta_j^{v}\lambda^{v}, \ v \in \{d,r\} %\label{slope_function}
%\end{align}
%where the parameter $\theta_j^{v}\in \mathbb{R}$ explicitly indicate willingness of generator $j$ to produce $g_j^{v}$ per unit price $\lambda^{v}$.

\subsection{Standard Two-stage Market}
\textcolor{black}{The role of participants in a standard market without any mitigation policy is studied extensively in the literature~\cite{Li_Na_linear_SFE,Joh2004Efficiency,baldick_linear_supply,pyou_discovery}. Here, we cite the results from~\cite{pyou_discovery} that analyze the role of strategic generators and inelastic demand in a standard two-stage market and use them as a benchmark to analyze the impact of a mitigation policy in the market.}

\subsubsection{Price-taking Participation and Competitive Equilibrium} 

For the individual incentive problem in a two-stage market, substituting the supply function~\eqref{gen_da_bid},\eqref{gen_rt_bid} in~\eqref{generator_profit}, we get 
\begin{align}
    \!\!\pi_j(\theta_j^{d},\theta_j^{r};\!\lambda^{d},\lambda^{r})\! = & \ \! \theta_j^{d}{\lambda^{d}}^2 \!\!+\! \theta_j^{r}{\lambda^{r}}^2\!\!-\!  \frac{c_j}{2}\!\!\left(\theta_j^{d}{\lambda^{d}} \!+\! \theta_j^{r}{\lambda^{r}}\right)^2\label{generator_price_taking_profit_bids_wout_mpm}
\end{align}
and the individual problem for price-taking generator $j$ is:
\begin{align}\label{generator_price_taking_profit_final_wout_mpm}
    &\max_{\theta_j^{d},\theta_j^{r}} \  \pi_j(\theta_j^{d},\theta_j^{r};\!\lambda^{d},\lambda^{r})
\end{align}
The individual problem for load $l$ is given by~\eqref{load_payment_bids}. Given the prices, $\lambda^d,\lambda^r$, we next characterize the resulting competitive equilibrium due to competition between price-taking participants.
\begin{theorem}[Proposition  1~\cite{pyou_discovery}] \label{comp_eqbm_wout_mpm} A competitive equilibrium in a two-stage market exists and is explicitly given by
\begin{subequations}
\begin{align}
    &{\theta}_j^d+{\theta}_j^r = c_j^{-1}, \ {\theta}_j^d \ge 0, \ {\theta}_j^r \ge 0, \forall j \in \mathcal{G}\\
    & d_l^d+d_l^r = d_l, \forall l \in \mathcal{L}\\
    & \lambda^d = \lambda^r = \frac{d}{\sum_{j\in\mathcal{G}}c_j^{-1}}
\end{align}
\end{subequations}
\end{theorem}
The resulting competitive equilibrium solves the social planner problem~\eqref{planner_problem}. Moreover, it exists non-uniquely, and there is no incentive for a load to allocate demand in the day-ahead market due to equal prices in two stages. 

\subsubsection{Price-Anticipating Participation and Nash Equilibrium}
The individual problem of price-anticipating generator $j$ and price-anticipating load $l$ is given by~\eqref{generator_strategic_profit_total} and~\eqref{load_strategic_payment}, respectively. We next characterize the resulting Nash equilibrium in the market.
\begin{theorem}[Proposition  4~\cite{pyou_discovery}] \label{strat_eqbm_wout_mpm}Assuming strategic generators are homogeneous $(c_j := c, \ \forall j \in  \mathcal{G})$ and make identical bids $(\theta_j^v := \theta^v, \ \forall j \in  \mathcal{G},  v \in \{d,r\})$ at equilibrium. If there are at least three firms, i.e., $G \ge 3$, a Nash equilibrium
in a two-stage market exists. Further, this equilibrium is unique and explicitly given by
\begin{align}
    &\!\!\!{\theta}_j^d \!=\!\! \frac{L(G-1)+1}{L(G-1)}\frac{G-2}{G-1}\frac{1}{c}, \ \! {\theta}_j^{r} \!\!=\!\! \frac{1}{L+1}\frac{(G-2)^2}{(G-1)^2}\frac{1}{c} \label{strat_eqbm_wout_mpm_eq.a}\\
    & \!\!\! d_l^{d} = \frac{L(G-1)+1}{L(L+1)(G-1)}d, \ d_l^{r} \!\!= d_l - d_l^d \label{strat_eqbm_wout_mpm_eq.b}\\
    & \lambda^{d} =  \frac{L}{L+1}\frac{G-1}{G-2}\frac{c}{G}d, \ \lambda^{r} =  \frac{G-1}{G-2}\frac{c}{G}d\label{strat_eqbm_wout_mpm_eq.c}
\end{align}\label{strat_eqbm_wout_mpm_eq}
\end{theorem}
The resulting Nash equilibrium exists uniquely, where price-anticipating loads anticipate the actions of generators and allocate demand to exploit lower prices in the day-ahead market. Thus prices are different in two stages. Moreover, the net demand allocation in the day-ahead and real-time markets follows 
\begin{align}
    \!\sum\nolimits_{l \in \mathcal{L}}d_l^d = d^d \in (0.5d,d), \ \sum\nolimits_{l \in \mathcal{L}}d_l^r \!=\! d^r \in (0,0.5d)
\end{align}
\subsection{Real-Time MPM Policy}

In this section, we first discuss the modified market model, the individual incentives of participants, and then characterize market equilibrium for a real-time MPM policy.

\subsubsection{Modeling Real-Time MPM Policy}
\begin{figure}
    \centering
    \includegraphics[width=0.8\linewidth]{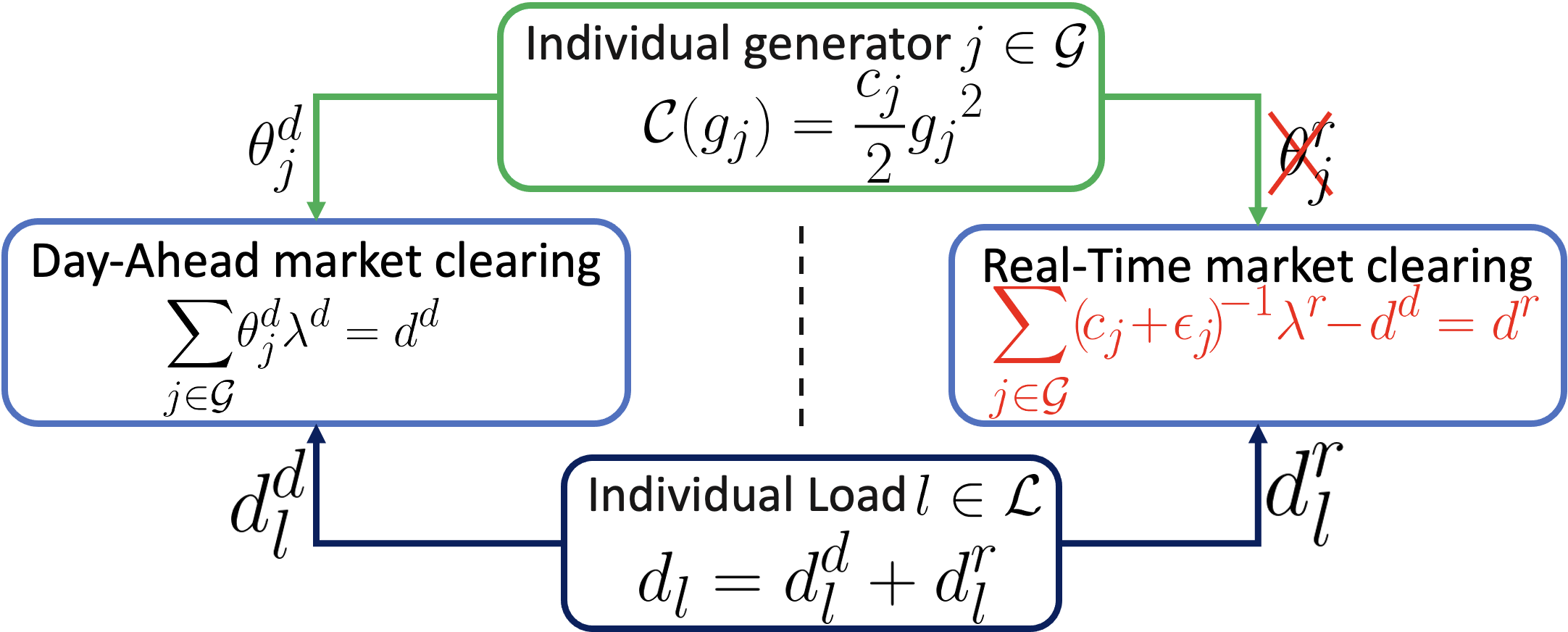}
    \caption{Two-stage Market Mechanism with Real-Time MPM}
    \label{fig:rt_mpm}
\end{figure}

In the case of a real-time MPM policy, the market ignores generators' bids in real-time, as shown in Figure~\ref{fig:rt_mpm}, and \textcolor{black}{roughly estimates the cost of dispatching generator $j$ with an error $\epsilon_j\ge0$}, given the day-ahead dispatch $g_j^{d}$ 
\textcolor{black}{\begin{align} \label{rt_true_dispatch}
        g_j^{r} = (c_j+\epsilon_j)^{-1}\lambda^{r}- g_j^{d}
\end{align}}
Using the two-stage generation and supply-demand balance~\eqref{two_stage_power_bal} and real-time dispatch~\eqref{rt_true_dispatch} we get
\textcolor{black}{\begin{align}\label{rt_true_prc}
    {\lambda^{r}} = \frac{d}{\sum\limits_{j\in\mathcal{G}}(c_j+\epsilon_j)^{-1}}
\end{align}}

\subsubsection{Price-taking Participation and Competitive Equilibrium} 

For the individual incentive problem in a two-stage market with real-time MPM policy, substituting the day-ahead supply function~\eqref{gen_da_bid}, real-time true dispatch condition~\eqref{rt_true_dispatch} and real-time clearing prices~\eqref{rt_true_prc} in~\eqref{generator_profit}, we get 
\textcolor{black}{\begin{align}
    \!\!\pi_j(\theta_j^{d},\!\lambda^{d})\! = & \ \! \theta_j^{d}{{\lambda\!}^{d}}^2 \!\!+\! \frac{d}{\sum\limits_{k\in\mathcal{G}}\!\!c_k^{-1}}\!\!\left(\!\omega_jd\!- \theta_j^{d}{{\lambda\!}^{d}}\!\right)\!\!-\!  \frac{c_j}{2}\!{\left(\omega_jd\right)}^2 \label{generator_price_taking_profit_bids_rt_mpm}
    %= &  \left(\!\!\lambda^{d} -\frac{d}{\sum\limits_{k\in\mathcal{G}}\!\!c_k^{-1}}\right)\!\theta_j^{d}{\lambda^{d}}+\!  \frac{c_j^{-1}}{2}\!\!{\left(\!\frac{d}{\sum\limits_{k\in\mathcal{G}}\!\!c_k^{-1}}\!\!\right)\!\!}^2
\end{align}}
where $\omega_j := \frac{(c_j+\epsilon_j)^{-1}}{\sum\nolimits_{k\in\mathcal{G}}(c_k+\epsilon_k)^{-1}}$. 
Hence, an individual problem of a price-taking generator is:
\begin{align}\label{generator_price_taking_profit_final_rt_mpm}
    &\max_{\theta_j^{d}} \  \pi_j(\theta_j^{d};\lambda^{d}) 
\end{align}
Similarly, substituting the clearing price~\eqref{rt_true_prc} in~\eqref{load_payment_bids} we get,
\textcolor{black}{\begin{align}\label{load_payment_intermediate_rt_mpm}
    & \rho_l(d_l^{d},\lambda^{d}) :=  \lambda^{d}d_l^{d} + \frac{d}{\sum_{k\in\mathcal{G}}(c_k+\epsilon_k)^{-1}}(d_l-d_l^{d})
\end{align}}
such that the individual problem for load $l$ is given by:
\begin{align}\label{load_payment_final_rt_mpm}
    \min_{d_l^{d}} \  &  \rho_l(d_l^{d};\lambda^{d})
\end{align}

%Given the market clearing prices $\lambda^{d}$ each generator $j \in \mathcal{G}$ maximize its profit~\eqref{generator_price_taking_profit_final_rt_mpm}. Similarly, given both clearing prices $\lambda^{d},\lambda^{r}$, each load $l\in\mathcal{L}$ minimizes its individual payment~\eqref{load_payment_final_rt_mpm}. 
The competition between price-taking participants for individual incentives leads to a set of competitive equilibria, as characterized below.

\textcolor{black}{\begin{theorem} \label{comp_eqbm_rt_mpm_Thrm1}
The competitive equilibrium in a two-stage market with a real-time MPM policy exists, and given by:
\begin{subequations}
\begin{align}
    & g_j^{d} + g_j^{r} \!= \frac{(c_j+\epsilon_j)^{-1}}{\sum\nolimits_{k\in\mathcal{G}}(c_k+\epsilon_k)^{-1}}d,  \ \theta_j^d \in \mathbb{R}_{\ge 0} \ \forall j\in\mathcal{G}\\
    & d_l^{d} + d_l^{r} = d_l, \ \!\! \forall l\in\mathcal{L}\\
    & \lambda^{d} = \lambda^{r} = \frac{d}{\sum\nolimits_{k\in\mathcal{G}}(c_k+\epsilon_k)^{-1}}
\end{align}
\label{competitive_eqbm_traditional_rt_mpm}
\end{subequations}
\end{theorem}}
\ifthenelse{\boolean{arxiv}}{We provide proof of the theorem in Appendix~\ref{appendix_comp_eqbm_rt_mpm}}{The proof is provided in~\cite{bansal2023market}}. At the competitive equilibrium, the market clearing prices are equal in the two stages, meaning there is no incentive for a load to allocate demand in the day-ahead market, e.g., current market practice. \textcolor{black}{However, the resulting equilibrium in Theorem~\ref{comp_eqbm_rt_mpm_Thrm1} is inefficient and does not always align with the social planner problem.}
\textcolor{black}{\begin{corollary}
    The competitive equilibrium in a two-stage market with a real-time MPM policy~\eqref{competitive_eqbm_traditional_rt_mpm} also solves  the social planner problem~\eqref{planner_problem} only when $\epsilon_j = 0 , \ \forall j \in \mathcal{G}$,
\end{corollary}}

\subsubsection{Price-Anticipating Participation and Nash Equilibrium}
The individual problem of each price-anticipating generator $j$, given by:
%\begin{subequations}\label{generator_strategic_profit_total_rt_mpm}
\begin{align}\label{generator_strategic_profit_rt_mpm}
    & \!\!\max_{\theta_j^{d},\lambda^d} \ \! \pi_j \! \left(\theta_j^{d},\lambda^{d}\!\left(\theta_j^{d};\overline{\theta}_{-j}^{d},d^{d}\right)\right) \ \textrm{ s.t. } \eqref{da_power_bal}
\end{align}
%\end{subequations}
where generator $j$ maximizes its profit in the two-stage market. The individual problem of price-anticipating load is:
%\begin{subequations}
\begin{align}\label{load_strategic_payment_rt_mpm}
    & \!\!\!\min_{d_l^{d},\lambda^d} \ \rho_l\!\left(d_l^{d}, \lambda^{d}\!\left(d_l^d;\theta_j^{d},\overline{d}_{-l}^{d}\right)\right) \ \textrm{ s.t. } \eqref{da_power_bal}
\end{align}
%\end{subequations}
where the load $l$ minimizes its payment in the market. 

We study the resulting sequential game where players anticipate each other actions and prices in the market, and the day-ahead clears before the real-time market. To this end, we analyze the game backward, starting from the real-time market, where prices are fixed due to MPM policy~\eqref{rt_true_prc}, followed by the day-ahead market, where participants make decisions for optimal individual incentives and compute the equilibrium path. Generators do not bid in real-time, but loads are allowed to bid in the market. However, load makes decisions simultaneously in the day-ahead market due to inelasticity, fixing their bids in the real-time market, which affects the two-stage market clearing. The following theorem characterizes the two-stage Nash equilibrium that satisfies the Definition~\eqref{market_eqbm}. 
\textcolor{black}{\begin{theorem} \label{strat_thrm_rt_mpm}
The Nash equilibrium in a two-stage market with a real-time MPM policy does not exist.
%\begin{subequations}
%\begin{align} \label{strat_eqbm_traditional_rt_mpm.b0}
%    & g_j^{d} = 0,  \ g_j^{r} = \frac{c_j^{-1}}{\sum_{j\in\mathcal{G}}c_j^{-1}}d, \ \theta_j^{d} = 0, \ \forall j\in\mathcal{G} \\ \label{strat_eqbm_traditional_rt_mpm.b1}
%    & d_l^{d} = 0, \ d_l^{r} = d_l, \ \forall l\in\mathcal{L} \\ \label{strat_eqbm_traditional_rt_mpm.b2}
%    & \lambda^{d} =  0, \lambda^{r} = \frac{1}{\sum_{j\in\mathcal{G}}c_j^{-1}}d   
%\end{align}\label{strat_eqbm_rt_mpm}
%\end{subequations}
\end{theorem}}
\ifthenelse{\boolean{arxiv}}{We provide proof of the theorem in Appendix~\ref{appendix_strat_eqbm_rt_mpm}}{We provide proof of the theorem in~\cite{bansal2023market}} and a brief insight below into the loss of equilibrium. The price-anticipating participants compete with each other to manipulate prices in the day-ahead given by~\eqref{da_power_bal}: 
\begin{align}
    \lambda^d = \frac{d^d}{\sum\limits_{j\in\mathcal{G}}\theta_j^d}  
\end{align}
while the prices in the real-time $\lambda^{r}$~\eqref{rt_true_prc} is fixed. \textcolor{black}{Loads bid decreasing quantities $d_l^{d}$ to reduce clearing prices in the day-ahead market and minimize the load payment. Simultaneously, generators bid decreasing parameter $\theta_j^{d}$ to increase clearing prices and maximize revenue. The competition between loads and generators for individual incentives in the day-ahead market drives all the demand to the real-time market, where generators operate truthfully. However, in our market mechanism, loads then have the incentive to deviate and allocate demand in the day ahead where prices are zero, meaning zero payment in the market, see Rule~\ref{assumpt_zero_prc}. Such unilateral load deviations result in deviations from generators to increase clearing prices in the day-ahead market. Therefore the equilibrium does not exist. Without such a market rule, the Nash equilibrium does exist with undefined clearing prices in the day-ahead and all demand allocated to the real-time market. Nevertheless, since day-ahead accounts for a majority of energy trades, the resulting equilibrium is undesirable.
}

%The competition between loads and generators for individual incentives in the day-ahead market drives all the demand to the real-time market, where generators operate truthfully. However, in our market mechanism, loads then have the incentive to deviate and allocate demand in the day ahead where prices are zero, meaning zero payment in the market, see Rule~\ref{assumpt_zero_prc}. Therefore the equilibrium does not exist. Without such market rules, the equilibrium does exist with undefined clearing prices in the day-ahead and all demand in the real-time. Since day-ahead accounts for majority of energy trades, we get rid of such undesirabel equilibriums in this study.  
%Although the Nash equilibrium~\eqref{strat_eqbm_rt_mpm} solves the social planner problem~\eqref{planner_problem}, the equilibrium is  undesirable from a market perspective since the current market practice shows majority of energy trades in the day-ahead market.

\subsection{Day-Ahead MPM Policy}

In this section, we define the individual incentive of participants and characterize market equilibrium for a day-ahead MPM policy.

\subsubsection{Modeling Day-Ahead MPM Policy}
In the case of a day-ahead MPM policy, as shown in Figure~\ref{fig:da_mpm}, the market ignores the generators' bids and \textcolor{black}{roughly estimates the cost of dispatching generator $j$ in the day-ahead with an error $\epsilon_j\ge0$}, as given by: 
\begin{figure}
    \centering
    \includegraphics[width=0.8\linewidth]{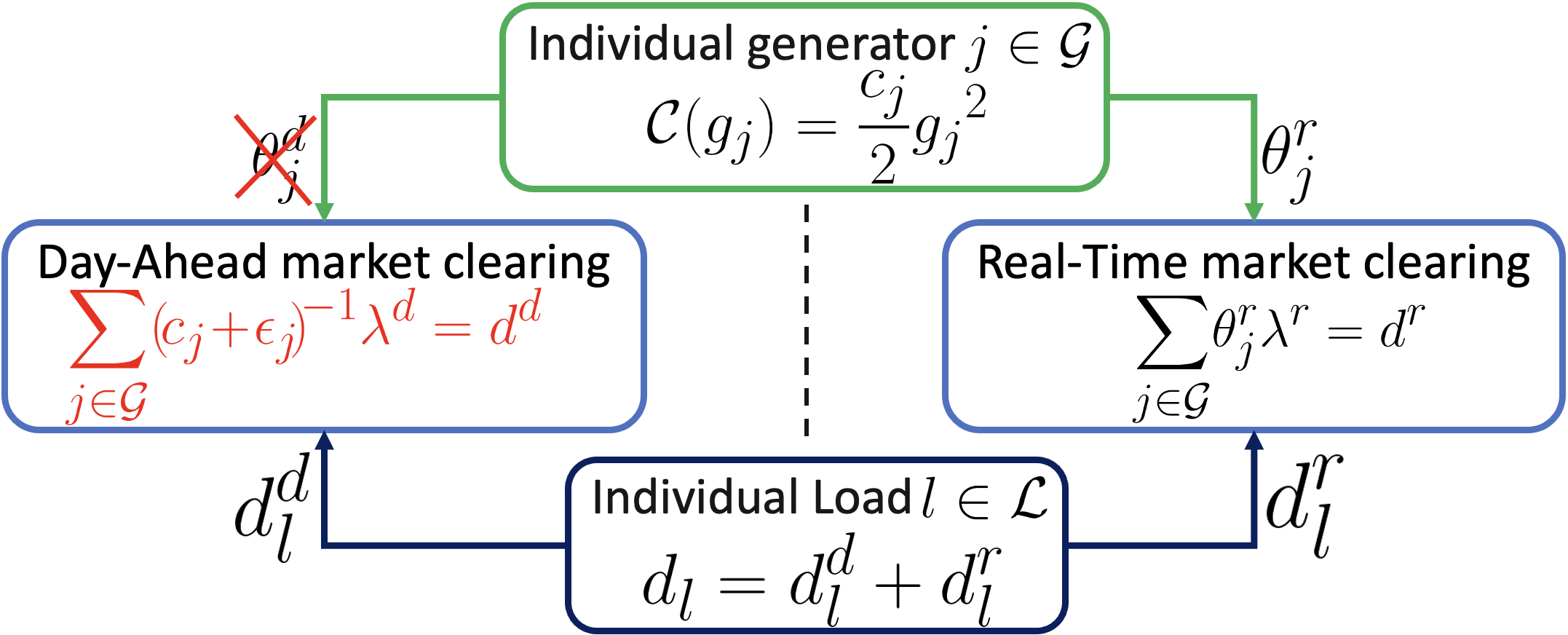}
    \caption{Two-stage Market Mechanism with Day-Ahead MPM}
    \label{fig:da_mpm}
\end{figure}
\textcolor{black}{\begin{align} \label{da_true_dispatch}
    g_j^{d} = (c_j+\epsilon_j)^{-1}\lambda^d
\end{align}}
Moreover, using day-ahead power balance constraint, we get
\textcolor{black}{\begin{align}\label{da_true_prc}
    \lambda^{d} = \frac{d^d}{\sum_{j\in\mathcal{G}}(c_j+\epsilon_j)^{-1}}
\end{align}}

\subsubsection{Price-taking Participation and Competitive Equilibrium}
For the individual incentive problem in a two-stage market with a day-ahead MPM policy, substituting the clearing price~\eqref{da_true_prc} in~\eqref{generator_profit}, we get 
\textcolor{black}{\begin{align}
    \! \pi_j(\theta_j^{r};\!\lambda^{r}) \!= \!\!(c_j\!+\!\epsilon_j)(\omega_j{d^d})^2\! +\! \theta_j^r{\lambda^{r}}^2 \!-\!  \frac{c_j}{2}\!\!\left(\!\omega_jd^d\! +\! \theta_j^r{\lambda^{r}}\!\right)^2 \label{generator_price_taking_profit_bids_da_mpm}
\end{align}}
where $\omega_j := \frac{(c_j+\epsilon_j)^{-1}}{\sum\nolimits_{k\in\mathcal{G}}(c_k+\epsilon_k)^{-1}}$. The individual problem for price-taking generator $j$ is:
\begin{align}\label{generator_price_taking_profit_final_da_mpm}
    &\max_{\theta_j^{r}} \  \pi_j(\theta_j^{r};\lambda^{r}) 
\end{align}
and the individual problem for load $l$ is given by~\eqref{load_payment_bids}.
%\begin{align}\label{load_payment_final_da_mpm}
%    \!\!\! \min_{d_l^{d}} \rho_l(d_l^{d};\lambda^{r}) \!=\! %\min_{d_l^{d}} \frac{d^d}{\sum_{j\in\mathcal{G}}(c_j\!+%\!\epsilon_j)^{-1}}d_l^{d} \!+\! \lambda^{r}(d_l-d_l^{d})
%\end{align}
%Each price-taking generator $j$ solves individual problem~\eqref{generator_price_taking_profit_final_da_mpm} and load $l$ solve individual problem~\eqref{load_payment_bids}. 
The resulting competitive equilibrium given the clearing prices $\lambda^{d}$ and $\lambda^{r}$ is characterized below.
\textcolor{black}{
\begin{theorem}\label{comp_eqbm_da_mpm_thrm}
The competitive equilibrium in a two-stage market with a day-ahead MPM policy exists and is given by:
\begin{subequations}
\begin{align}
    &\!\!\! g_j^{d} = \frac{(c_j+\epsilon_j)^{-1}}{\sum_{j\in\mathcal{G}}c_j^{-1}}d,\! \ g_j^{r} = \frac{\epsilon_jc_j^{-1}}{c_j\!+\!\epsilon_j}\frac{d}{\sum_{j\in\mathcal{G}}c_j^{-1}} \ \forall j\in\mathcal{G} \\
    &\!\!\! d_l^{d}\!+\! d_l^{r} \!=\! d_l, \ \!\forall l\in\mathcal{L},\! \  d^{d} \!=\! \frac{\sum\limits_{j\in\mathcal{G}}\!(c_j\!+\!\epsilon_j)^{-1}}{\sum\limits_{j\in\mathcal{G}}c_j^{-1}}d,\!\ d^{r} \!=\! d\!-\!d^d \label{comp_eqbm_da_mpm.b}\\
    &\!\! \theta_j^r = \frac{\epsilon_jc_j^{-1}}{c_j\!+\!\epsilon_j}, \ \lambda^{d} = \lambda^{r} = \frac{1}{\sum_{j\in\mathcal{G}}c_j^{-1}}d
\end{align}
\label{comp_eqbm_da_mpm}
\end{subequations}
\end{theorem}}
\ifthenelse{\boolean{arxiv}}{We provide proof of the theorem in Appendix~\ref{appendix_comp_eqbm_da_mpm}}{We provide proof of the theorem in~\cite{bansal2023market}}. Unlike the competitive equilibrium for a real-time MPM policy in~\eqref{competitive_eqbm_traditional_rt_mpm} with equal prices across stages, the loads at equilibrium~\eqref{comp_eqbm_da_mpm} \textcolor{black}{allocate a majority of the demand} in the day-ahead. The incentive for day-ahead demand allocation is a desired market outcome and is not generally satisfied by other market mechanisms. %Although $d^r=0$ and $\sum_j\theta_j^r = 0$, unlike in Theorem~\ref{strat_thrm_rt_mpm}, loads do not deviate from the equilibrium by shifting demand to the real-time market. 
The resulting equilibrium exists as price-taking loads do not anticipate the effect of their bid on the market prices, meaning the payment remains the same for any allocation across the two stages. \textcolor{black}{Moreover, the market outcome~\eqref{comp_eqbm_da_mpm} solves the social planner problem~\eqref{planner_problem}.}

\subsubsection{Price-Anticipating Participation and Nash Equilibrium}

The individual problem of each price-anticipating generator $j$, given by:
%\begin{subequations}%%\label{generator_strategic_profit_total_da_mpm}
\begin{align}\label{generator_strategic_profit_da_mpm}
    & \!\!\max_{\theta_j^{r},\lambda^r} \ \! \pi_j \! \left(\theta_j^{r},\lambda^{d}\!\left(d^{d}\right),\lambda^{r}\!\!\left(\theta_j^{r};\overline{\theta}_{-j}^{r},d^{r}\right)\right) \textrm{ s.t. } \eqref{rt_power_bal}
\end{align}
%\end{subequations}
where generator $j$ maximizes its profit in the market. The individual problem of price-anticipating load $l$, is given by:
%\begin{subequations}
\begin{align}\label{load_strategic_payment_da_mpm}
    & \!\!\!\min_{d_l^{d},\lambda^r} \ \rho_l\!\left(d_l^{d}, \lambda^{d}\!\left(d_l^{d};\overline{d}_{-l}^d\right),\lambda^{r}\!\!\left(d_l^d;\theta_j^{r},\overline{d}_{-l}^{r}\right)\right) \textrm{ s.t. } \eqref{rt_power_bal}
\end{align}
%\end{subequations}
where load $l$ minimizes its payment in the market. 

In the market model with a day-ahead MPM policy, generators make decisions in real-time while load can make decisions in the day-ahead. The resulting two-stage sequential game is essentially a leader-follower Stackelberg-Nash game, where generators are followers in the real-time market and loads are leaders in the day-ahead market, and each participant in their respective groups competes amongst themselves in a Nash game. We follow the terminology used in~\cite{stack_nash_tnn} to describe similar formulations in different markets. For the closed form solution, we assume that generators are homogeneous in the sense that they share the same cost coefficient, i.e. $c_j =: c, \ \forall j \in \mathcal{G}$ and bid symmetrically in the market, i.e. $\theta_j^{r} =: \theta^{r}, \ \forall j \in \mathcal{G}$. Under these assumptions, the Nash equilibrium is characterized below. 

\textcolor{black}{
\begin{theorem}\label{strat_thrm_da_mpm}
Assume that generators are homogeneous and bid symmetrically in the market. Also, assume that estimation error is same for homogeneous generators, i.e. $\epsilon_j:=\epsilon, \ \forall j \in \mathcal{G}$. If more than two generators are participating in the market i.e., $G \ge 3$ and the number of individual loads participating in the market satisfies $\frac{1}{L} > \frac{c-\epsilon(G-2)}{(c+\epsilon)(G-2)}$, then the symmetric Nash equilibrium in a two-stage market with a day-ahead MPM policy exists uniquely as:
\begin{subequations}
\begin{align}
    & \!\!\! g_j^{d} \!=\! \frac{c}{c\!+\!\epsilon}\frac{L}{L\!+\!1}\frac{G\!-\!1}{G\!-\!2}\frac{d}{G},\! \ g_j^{r} \!\!=\!\! \left(\!1\!-\!\frac{c}{c\!+\!\epsilon}\frac{L}{L\!+\!1}\frac{G\!-\!1}{G\!-\!2}\!\right)\frac{d}{G} \label{strat_eqbm_da_mpm.a}\\
    & \!\!\! d_l^{d} \!=\! \frac{c}{c\!+\!\epsilon}\frac{1}{L\!+\!1}\frac{G\!-\!1}{G\!-\!2}d,\! \ d_l^{r} \!\!= \!\!\left(\!\!d_l \!-\! \frac{c}{c\!+\!\epsilon}\frac{1}{L\!+\!1}\frac{G\!-\!1}{G\!-\!2}d\right) \label{strat_eqbm_da_mpm.b}\\
    & \theta^{r} \!=\! \frac{1}{c}\left(\frac{G\!-\!2}{G\!-\!1}\!-\!\frac{c}{c\!+\!\epsilon}\frac{L}{L\!+\!1}\right) \label{strat_eqbm_da_mpm.c}\\
    & \lambda^{d} =  \frac{L}{L\!+\!1}\frac{G\!-\!1}{G\!-\!2}\frac{c}{G}d, \ \lambda^{r} =  \frac{G\!-\!1}{G\!-\!2}\frac{c}{G}d.  \label{strat_eqbm_da_mpm.d} 
\end{align}\label{strat_eqbm_da_mpm}%
\end{subequations}%
Moreover, for $\frac{1}{L} \le \frac{c-\epsilon(G-2)}{(c+\epsilon)(G-2)}$, a symmetric equilibrium does not exist. 
\end{theorem}}

\ifthenelse{\boolean{arxiv}}{We provide proof of the Theorem in Appendix~\ref{appendix_strat_eqbm_da_mpm}}{We provide proof of the theorem in~\cite{bansal2023market}}. Unlike the market with a real-time MPM policy, the Nash equilibrium exists in the market with a day-ahead MPM policy. However, it requires restrictive conditions on the number of participants in the market and may not even exist in other cases. We discuss these cases with no symmetric Nash equilibrium and provide intuition into participants' behavior in the market: 

\begin{boxlabel}
    \item $\mathbf{\frac{1}{L} < \frac{c-\epsilon(G-2)}{(c+\epsilon)(G-2)}}:$ In this case, the net demand is negative in the real-time market. The first order condition implies that each generator $j$ acts as load, paying $\lambda^{r}g_j^{r}$ as part of the market settlement since their optimal bid $\theta_j^r<0$ and the real-time clearing price $\lambda^{r}>0$. However, if the generators bid $\theta_j^r>0$, then the linear supply function implies that each generator $j$ dispatch $g_j^{r}<0$ at the clearing prices $\lambda^{r} <0$ earning revenue in the market. However, this is not desirable from a load perspective since they are making payments in the market and they have the incentive to deviate to minimize their payment. Hence, symmetric equilibrium with negative demand in the real-time market does not exist as the symmetric bid $\theta_j^r > 0$ does not satisfy the first-order condition. The dependence of the individual bid $\theta_j^r$ on the given bids from other participants makes the closed-form analysis challenging, and any guarantee of the existence of equilibrium is hard.

    \item {$\mathbf{\frac{1}{L} = \frac{c-\epsilon(G-2)}{(c+\epsilon)(G-2)}}:$} In this case, no symmetric Nash equilibrium exists. %and the net demand is zero in the real-time market. 
Loads take advantage of the truthful participation of generators in day-ahead market and their ability to anticipate impact of bids on the clearing prices. Regardless of generators' bids, loads have the incentive to deviate by allocating demand in the real-time market with a lower clearing price.
\end{boxlabel}

\begin{corollary}
\textcolor{black}{For ${\frac{1}{L} > \frac{c-\epsilon(G-2)}{(c+\epsilon)(G-2)}}$, at the Nash equilibrium~\eqref{strat_eqbm_da_mpm} in a two-stage market with a day-ahead MPM policy, the demand allocation is given by:
    \begin{subequations}
    \begin{eqnarray}
        & \sum_{l\in\mathcal{L}} d_l^{d} = d^d = \frac{c}{c+\epsilon}\frac{L}{L+1}\frac{G-1}{G-2}d \\ 
\label{load_allocation_day_ahead}
    & \sum_{l\in\mathcal{L}}  d_l^{r} = d^r = \left(1- \frac{c}{c+\epsilon}\frac{L}{L+1}\frac{G-1}{G-2}d\right)d 
\label{load_allocation_real_time}
    \end{eqnarray}
    \end{subequations}
}
Assuming $\epsilon = 0$, the following relation holds,
\[
    d^d \in (0.5d,d), \ d^r \in (0,0.5d)
\]
\end{corollary}

%The proof uses the fact that $L<G-2$ implies that $d^{r}>0$.
\section{Equilibrium Analysis}\label{sec_4}

In this section, we study the properties of market equilibrium under the proposed policy framework and compare it with the standard market equilibrium.
\subsection{Comparison of a stage-wise MPM Policy}

\textcolor{black}{An MPM policy in real-time either results in an inefficient market outcome at the competitive equilibrium or leads to no Nash equilibrium. However, an MPM policy in the day-ahead leads to a stable market outcome that is robust to price manipulations, e.g. see Nash equilibrium~\eqref{strat_eqbm_da_mpm}. Despite errors in cost estimations, the competitive equilibrium is efficient~\eqref{comp_eqbm_da_mpm}.} This is summarized in Table~\ref{tab:eqbm_table}. 

We further analyze the case of a day-ahead MPM policy to study the strategic behavior of participants while regarding the respective competitive equilibrium in Theorem~\ref{comp_eqbm_da_mpm_thrm} as a benchmark. In the case of a day-ahead MPM policy, loads act as leaders in the day-ahead and generators as followers in real-time. The generator bids to manipulate prices leading to inflated prices in real-time~\eqref{strat_eqbm_da_mpm.d} while the load shifts its allocation in the day-ahead~\eqref{strat_eqbm_da_mpm.b}, increasing prices in the day-ahead market. Though the market equilibrium deviates from the competitive equilibrium~\eqref{comp_eqbm_da_mpm}, the social cost remains the same due to the homogeneous participation of generators. Table~\ref{tab:cost_profit_table} summarizes the aggregate profit and aggregate payment of generators and loads, respectively.
\begin{table}[!t]
  \caption{Competitive Equilibrium (CE) and Nash Equilibrium (NE) with a stage-wise MPM policy}
  \label{tab:eqbm_table}
  \centering
  \begin{tabular}{c|c|c}
    \hline
    Instance & Real-Time MPM & Day-Ahead MPM\\
    %\midrule
    \hline
    %\vspace{0.8 pt}
    \multirow{3}{*}{CE} & {Non-unique equilibrium} &{Unique equilibrium}\\
    & \textcolor{black}{Do not achieve social cost} & \textcolor{black}{Achieve social cost}\\
    & Arbitrary demand allocation  & Higher demand in day-ahead \\
    \hline
    %\noalign{\smallskip}
    %\vspace{0.8 pt}
    \multirow{3}{*}{NE} & Does not exist & Symmetric equilibrium\\
    & - & Social Cost same as CE\\
    & - & Extra constraints on players\\
    \hline
\end{tabular}
\end{table}
\textcolor{black}{
\begin{corollary}
    For $L< G-2$, the aggregate payment of loads and aggregate profit of generators at symmetric Nash equilibrium~\eqref{strat_eqbm_da_mpm} is less than that at respective competitive equilibrium~\eqref{comp_eqbm_da_mpm}. Moreover, for $L\ge G-2$ and ${\frac{1}{L} > \frac{c-\epsilon(G-2)}{(c+\epsilon)(G-2)}}$, the aggregate payment of loads and aggregate profit of generators at symmetric Nash equilibrium~\eqref{strat_eqbm_da_mpm} is greater than that at respective competitive equilibrium~\eqref{comp_eqbm_da_mpm}. 
\end{corollary}
}
The corollary follows from comparing the aggregate profit (payment) at Nash equilibrium to that at competitive equilibrium in Table~\ref{tab:cost_profit_table} for $L<G-2$.

\subsection{Equilibrium comparison with a standard market }

In this section, we compare the equilibrium in a day-ahead MPM policy market to a standard market. The social cost at the competitive equilibrium remains the same for the two markets with equal prices in the two stages. However, unlike in the case of a day-ahead MPM policy, the competitive equilibrium in Theorem~\ref{comp_eqbm_wout_mpm} exists non-uniquely and there is no incentive for a load to allocate demand in the day-ahead market. 

Interestingly, at Nash equilibrium prices in the two stages are the same for a day-ahead MPM policy market~\eqref{strat_eqbm_da_mpm.d} and a standard market~\eqref{strat_eqbm_wout_mpm_eq.c}. \textcolor{black}{Furthermore, an error in the estimation of the cost of dispatching generators does not impact market prices due to the participation of homogeneous generators.} However, the dispatch of generators and allocation of demand is different in the two market settings due to a leader-follower structure between participants in the market with a day-ahead MPM policy. %Assuming $L<G-2$, at Nash equilibrium we have 
%\[
%    d^d_{DA-MPM} -d^d_{Standard} =  \frac{1}{L+1}\left(\frac{L}{G-2} - \frac{1}{G-1}\right)d > 0
%\]
%meaning, a load is able to exploit the market further by allocating more demand in the day-ahead market with a lower clearing price. 
To understand the impact of price-anticipating participants on market equilibrium, we compare the aggregate profit (payment) in Table~\ref{tab:cost_profit_table} and \ref{tab:cost_profit_table_wout_mpm}, respectively.

\begin{table}
  \caption{\textcolor{black}{Comparison between Competitive Equilibrium (CE) and Nash Equilibrium (NE) in a market with a day-ahead MPM policy}}
  \label{tab:cost_profit_table}
  \centering
  \resizebox{\linewidth}{!}{\begin{tabular}{c|c|c}
    \hline
    Case & Generators total profit & Loads total payment\\
    %\midrule
    \hline
    \vspace{0.8 pt}
    CE \!\!\!&\!\!\!$\frac{1}{2}\frac{c}{G}d^2$&$\frac{c}{G}d^2$\\
    %\hline
    %\noalign{\smallskip}
    \vspace{0.8 pt}
    NE \!\!\!&\!\!\!$\frac{1}{2}\!\frac{c}{G}d^2\!\!\left(\!\frac{G}{G-2}\!-\!\frac{c}{c+\epsilon}\frac{(G-1)^2}{(G-2)^2}\frac{2L}{(L+1)^2}\!\!\right)$\!\!\!&\!\!\!$\frac{c}{G}d^2\!\!\left(\!\frac{G-1}{G-2}\!-\!\frac{c}{c+\epsilon}\frac{(G-1)^2}{(G-2)^2}\frac{L}{(L+1)^2}\!\!\right)$\!\!\!\!\!\\
    %\midrule
    \hline
\end{tabular}}
\end{table}

\begin{table}
  \caption{Comparison between Competitive Equilibrium (CE) and Nash Equilibrium (NE) in a standard market}
  \label{tab:cost_profit_table_wout_mpm}
  \centering
  \resizebox{\linewidth}{!}{\begin{tabular}{c|c|c}
    \hline
    Case & Generators total profit & Loads total payment\\
    %\midrule
    \hline
    \vspace{0.8 pt}
    CE \!\!\!&\!\!\! $\frac{1}{2}\frac{c}{G}d^2$&$\frac{c}{G}d^2$\\
    NE \!\!\!&\!\!\!$\frac{1}{2}\!\frac{c}{G}d^2\!\!\left(\!\frac{G}{G-2}\!-\!\frac{2L(G-1)+2}{(L+1)^2(G-2)}\!\right)$\!\!\!&\!\!\!$\frac{c}{G}d^2\!\!\left(\!\frac{G-1}{G-2}\!-\!\frac{L(G-1)+1}{(L+1)^2(G-2)}\!\right)$\!\!\!\!\!\\
    %\midrule
    \hline
\end{tabular}}
\end{table}

We restrict our comparison for $\frac{1}{L} > \frac{c-\epsilon(G-2)}{(c+\epsilon)(G-2)}$ only since the Nash equilibrium in Theorem~\ref{strat_thrm_da_mpm} does not exist otherwise. In particular, for $L = G-3$ the aggregate profit (payment) as shown in row 2 of Table~\ref{tab:cost_profit_table_wout_mpm} at the Nash equilibrium in Theorem~\ref{strat_eqbm_wout_mpm} equals to that of the competitive equilibrium. However, for $L < G-3$ the aggregate profit (payment) at Nash equilibrium is always less than the competitive equilibrium, meaning the loads are winners. The change in the normalized aggregate profit (payment) at the Nash equilibrium between a market with a day-ahead MPM policy and a standard market is given by 
\color{black}\[
    \frac{2}{(L+1)^2}\frac{1}{G-2} \left( 1 - \frac{L}{G-2}-L  + \frac{\epsilon}{c+\epsilon}L(G-1)\right) 
\]\color{black}
where profit (payment) is normalized with the competitive equilibrium. The difference depends on the number of participants and as the number of participants increases, the difference tends to $0$, since the Nash equilibrium in both cases approaches the competitive equilibrium, respectively.
\begin{figure*}
  \centering
  \includegraphics[width=1.\linewidth]{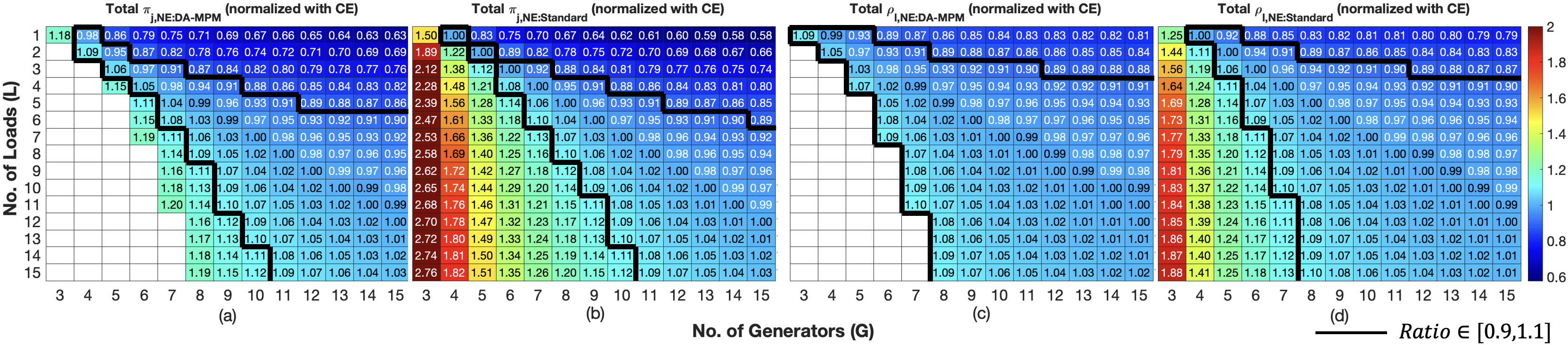}
  \caption{\textcolor{black}{Total profit and total payment at Nash Equilibrium (NE) normalized with competitive equilibrium (CE): total profit in (a) day-ahead MPM (DA-MPM) and (b) standard markets, and total payment in (c) day-ahead MPM (DA-MPM) and (d) standard markets; white cells denote no equilibrium}}
  \label{fig:ratio}
\end{figure*}

Figure~\ref{fig:ratio} compares the total profit (payment) normalized with competitive equilibrium for a day-ahead MPM (DA-MPM) policy market and a standard market for cost estimation error $\epsilon = 0.1$, respectively, as we change the number of loads ($l \in \mathcal{L}, \ L \in \{1,\ldots,G-3\}$), and generators ($j\in\mathcal{G},  \ G \in \{4,\ldots,20\}$). The ratio decreases monotonically as the number of generators increases, meaning the increased competition between more generators to meet the inelastic demand gives more power to loads, allowing them to reduce their payment even further, as shown by the horizontal rows in all panels in Figure~\ref{fig:ratio}. Furthermore, the ratio increases monotonically as the number of loads increases (for a large enough number of generators), meaning the market power shifts between loads and generators, as shown by the vertical color columns in panels (a) and (b) in Figure~\ref{fig:ratio}. In particular, in both markets, we observe a reversal in power, e.g., for a large number of loads generators make a higher profit at the expense of loads in the market and vice versa, as shown in panels (b) and (d) in the Figure~\ref{fig:ratio}. 

\textcolor{black}{Additionally, implementing the day-ahead MPM policy helps reduce market power. This leads to a total profit (payment) at Nash equilibrium that is closer to competitive equilibrium levels than what is observed in standard markets, as demonstrated in panels (a) and (c) for profit and panels (b) and (d) for payment in Figure~\ref{fig:ratio}.} Unfortunately, with a day-ahead MPM policy, the equilibrium does not always exists as shown by white-colored cells in panels (a) and (c). Finally, in the limit $L \rightarrow \infty \implies G \rightarrow \infty$, the Nash equilibrium converges to competitive equilibrium, also shown in Table~\ref{tab:cost_profit_table}. 

\section{Numerical Study}\label{sec_5}

%We now conduct a numerical case study to understand the impact of market participants on the Nash equilibrium in Theorem~\ref{strat_thrm_da_mpm}. We consider a two-stage market with a day-ahead market power mitigation policy, where multiple price-anticipating homogeneous generators and price-anticipating loads compete for individual incentives. For ease of comparison, we assume a fixed inelastic aggregate demand $d = 300 MW$. 
We now investigate how the cost estimation error, heterogeneity in cost coefficients, and load size affect individual incentives at Nash equilibrium in the market with a day-ahead MPM. We overcome the theoretical complexity of the closed-form analysis and run numerical best-response studies to understand the impact on market equilibrium. To this end, we consider the case of 2 price-anticipating loads and 5 price-anticipating generators in a two-stage market. 
\begin{figure}
  \centering
  \includegraphics[width=0.9\linewidth]{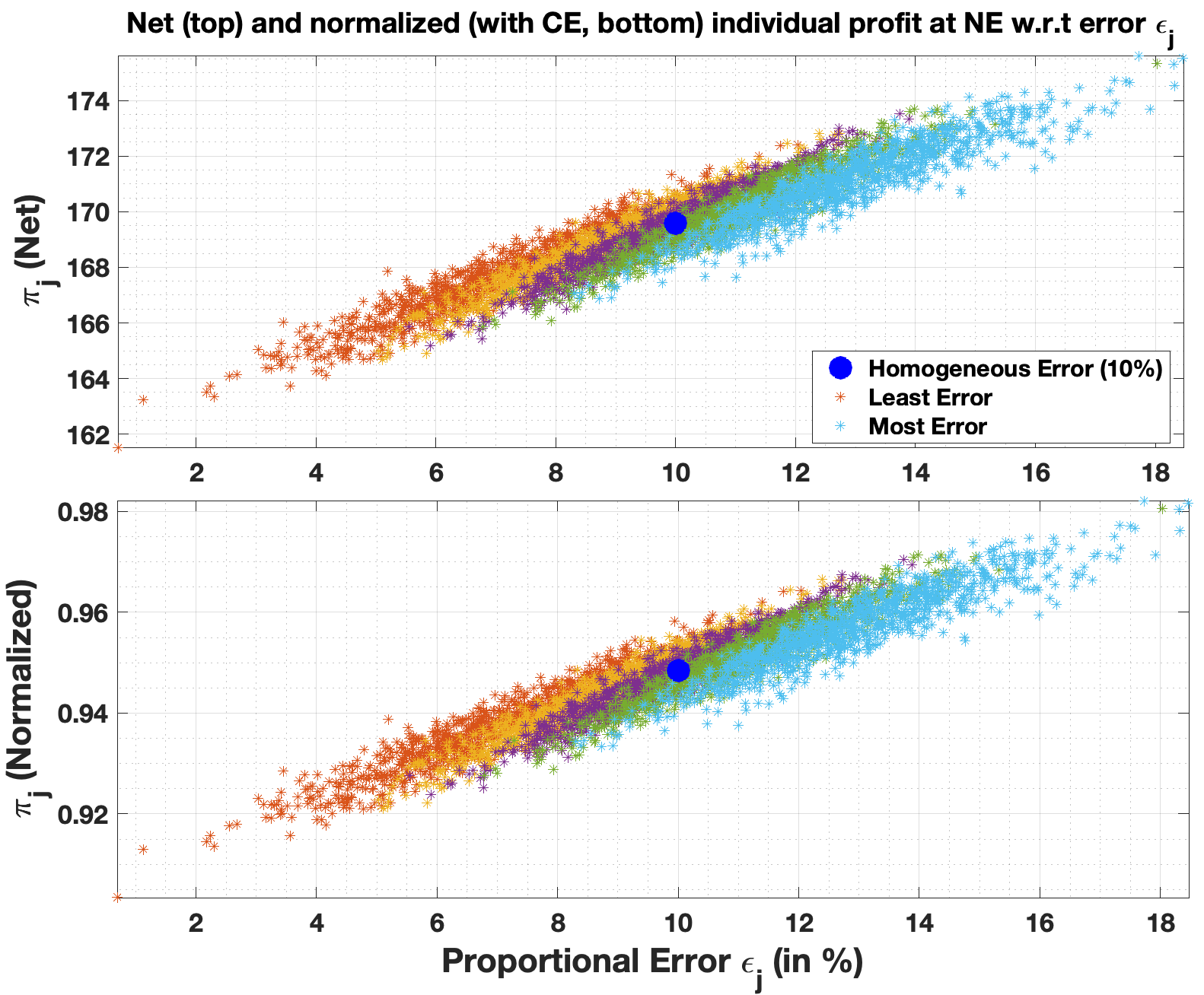}
  \caption{\textcolor{black}{Net (top) and normalized (bottom) individual profit at Nash Equilibrium (NE) normalized with competitive Equilibrium (CE) w.r.t proportional error $\epsilon_j$ in cost estimation of generators}}
  \label{fig:error_profit}
\end{figure}

\textcolor{black}{
The individual aggregate inelastic load is given by $d_l = [99.4, 199.6]^T \ MW$ from the Pennsylvania, New Jersey, and Maryland (PJM) data miner day-ahead demand bids~\cite{pjmdata}. For each generator $j$ with a truthful cost coefficient $c_j = 0.1 \$/MW^2,  \ \forall j \in \mathcal{G}$ corresponding to the cost coefficients from the IEEE 300-bus system~\cite{matpower}. We assume a proportional error $\epsilon_j = \delta_jc_j$ such that estimated cost coefficient is given by $\hat{c}_j = c_j(1+\delta_j), \ \forall j \in \mathcal{G}$. The cost estimation error of generators are sampled $10,000$ times from a Gaussian distribution with mean $10\%$ and variance $2.5\%$, i.e. $\delta_j \sim N(0.1,0.025) \ \forall j \in \{1,...,5\}$. The top and bottom panel in Figure~\ref{fig:error_profit} plots the net profit and the normalized profit (normalized with the competitive equilibrium) at Nash equilibrium, respectively. An increase in estimation error results in a higher net profit at Nash equilibrium, as shown in the  top panel in Figure~\ref{fig:error_profit}. Furthermore, errors in cost estimation also mitigate the market power of loads with profits closer to the competitive one.}
\begin{figure}
  \centering
  \includegraphics[width=0.9\linewidth]{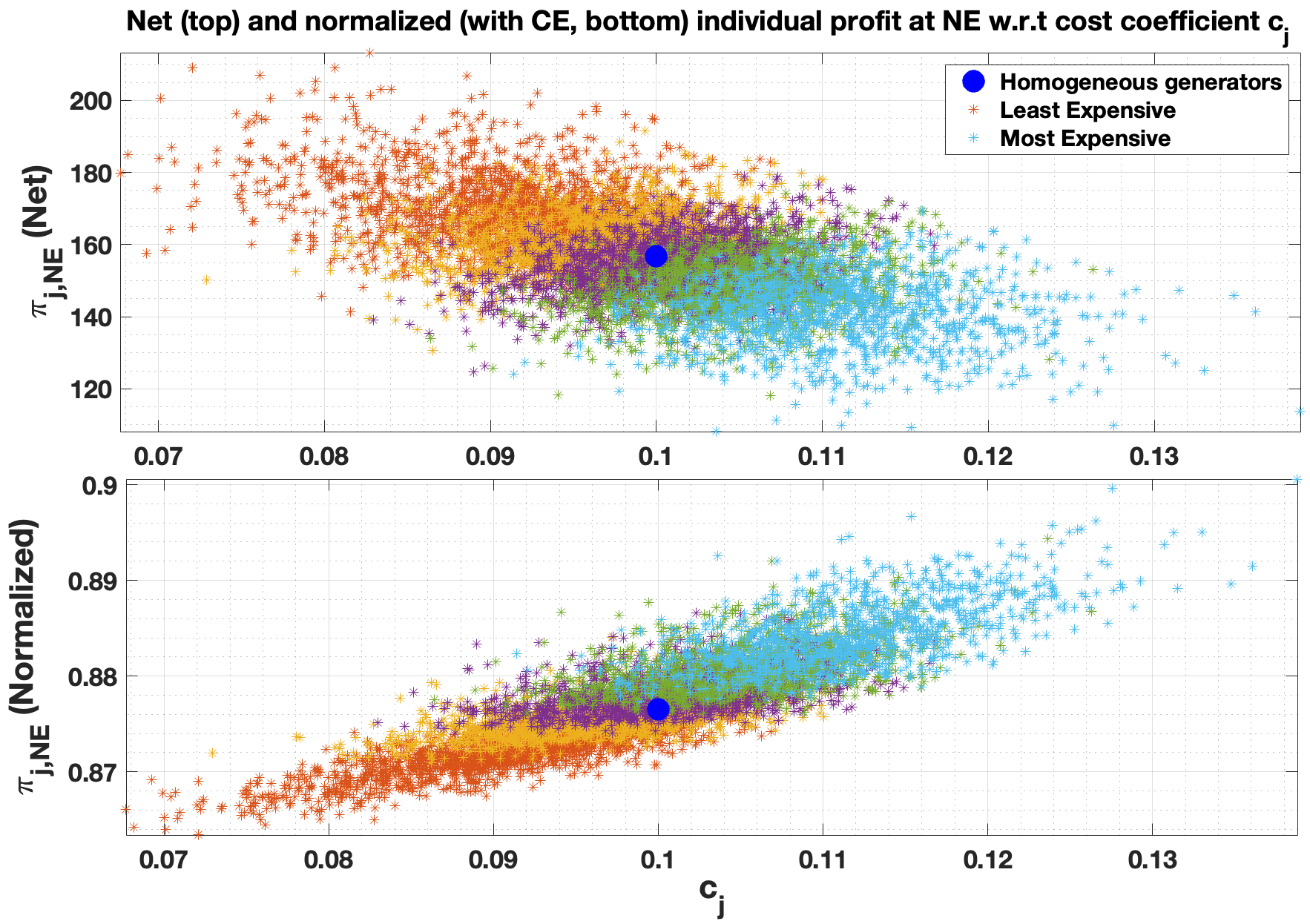}
  \caption{Net (top) and normalized (bottom) individual profit at Nash Equilibrium (NE) normalized with competitive Equilibrium (CE) w.r.t cost coefficient of generators for a DA-MPM policy}
  \label{fig:norm_profit}
\end{figure}

\textcolor{black}{
We next analyze the impact of heterogeneity in cost coefficients on market equilibrium. For ease of exposition, we assume that the cost estimation error $\epsilon_j = 0 \ \forall j \in \{1,...,5\}$. Our analysis is focused on capturing the qualitative impact of heterogeneity in cost coefficients on system-level market power. To this end, we choose a Gaussian distribution to model the uncertainty in the market operator's estimate for generators' truthful cost as the first step toward understanding the potential impact.} The cost coefficients of generators are sampled $10,000$ times from a Gaussian distribution with mean $0.1$ and sample variance $0.001$ for a sample of cost coefficients from the IEEE 300-bus system~\cite{matpower}, i.e. $c_j \sim N(0.1,0.001), \ \forall j \in \{1,...,5\}$. The top and bottom panel in Figure~\ref{fig:norm_profit} plots the absolute profit and the normalized profit (normalized with the competitive equilibrium) at Nash equilibrium, respectively. The cheaper generators earn a higher profit when compared with the expensive generators with higher cost coefficients at Nash equilibrium. However, the normalized profit ratio in the bottom panel shows that expensive generators have a higher value than cheaper ones, meaning that though expensive generators have lower absolute profit, these are the least exploited in the market. \textcolor{black}{We hypothesize that such a non-trivial behavior is related to the nature of competition between strategic generators instead of an effect of a day-ahead MPM policy. We admit that a closed-form analysis is theoretically complex, and we do not have a thorough mechanism to validate our hypothesis.} %This result is counter-intuitive.
\begin{figure}
  \centering
  \includegraphics[width=0.85\linewidth]{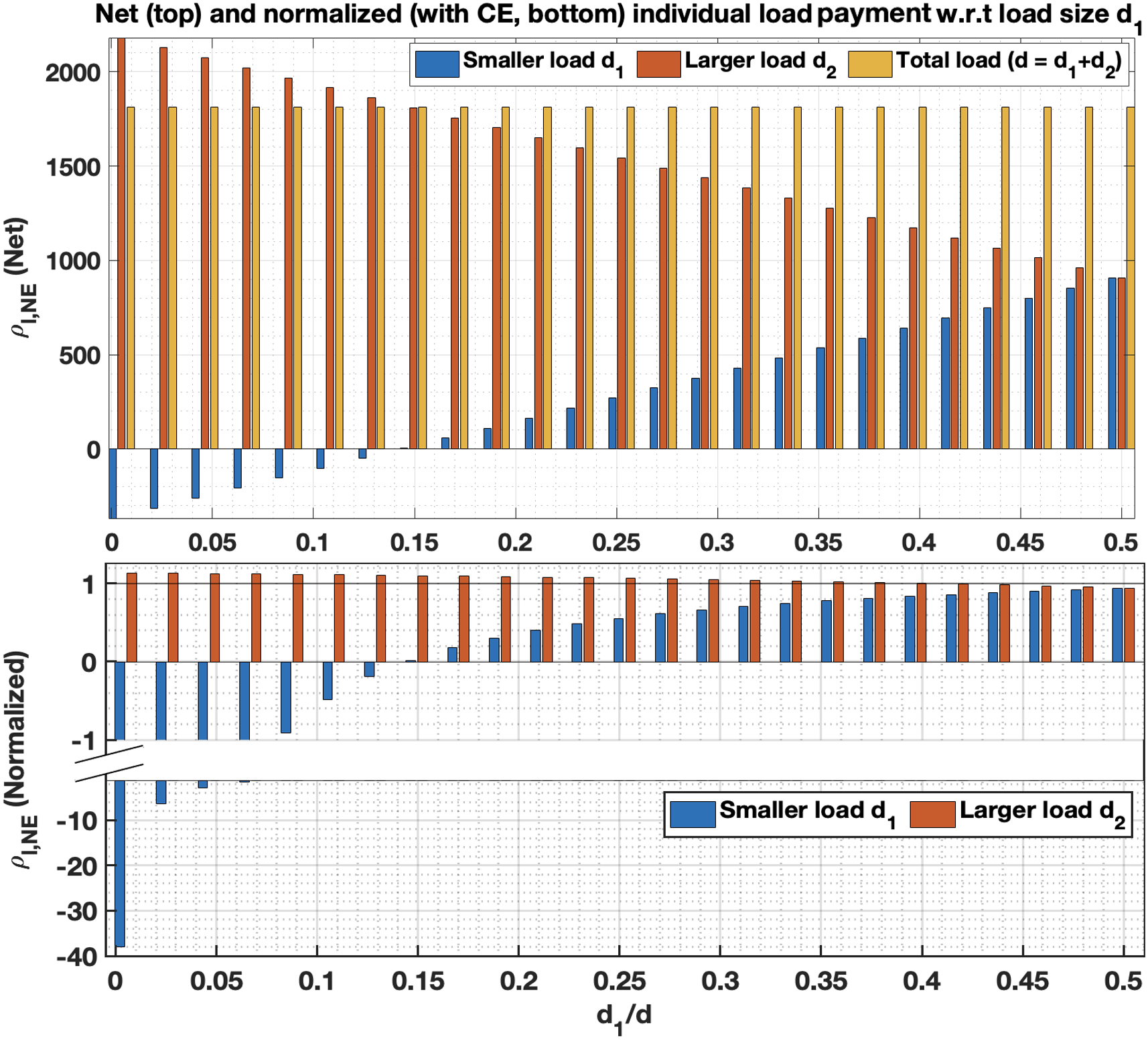}
  \caption{Net (top) and normalized (bottom) load individual payment (bottom) at Nash Equilibrium (NE) normalized with competitive Equilibrium (CE) w.r.t size of smaller load $d_1$, $d_1<d_2, \ d_1+d_2 = d$, for a DA-MPM policy}
  \label{fig:abs_payment}
\end{figure}

In Figure~\ref{fig:abs_payment} we show the absolute (top panel) and normalized (bottom panel) load payment w.r.t smaller load size. For this, we keep the same number of loads and generators in the market with varying load sizes for fixed net demand. We again sample cost coefficients from Gaussian distribution with mean $0.1$ and sample variance $0.001$ for a sample of cost coefficients from the IEEE 300-bus system~\cite{matpower}, i.e. $c_j \sim N(0.1,0.001), \ \forall j \in \{1,...,5\}$. \textcolor{black}{The cost estimation error $\epsilon_j = 0, \ \forall j \in \{1,...,5\}$.} The top panel shows that though the net load payment remains the same as we change the size of the load, the smaller load may even make a profit in the market at the expense of a higher load. More formally to develop intuition, in the case of homogeneous generators, the normalized payment ratio for individual load at Nash equilibrium in Theorem~\ref{strat_thrm_da_mpm}, is given by 
\[
    \frac{G-1}{G-2}\left(1-\frac{1}{(L+1)^2}\frac{G-1}{G-2}\frac{d}{d_l}\right)
\]
which is negative for a sufficiently small load. In particular, the smaller load has a negative normalized ratio at the expense of a higher load (a ratio greater than 1), as shown in the bottom panel of Figure~\ref{fig:abs_payment}. The larger load makes more payment at Nash equilibrium than at the competitive equilibrium, while the aggregate payment of the set of loads is still less than at the competitive equilibrium. Though the heterogeneity in load size does not affect the net payment or the group behavior in the market, a smaller load makes negative payments at the expense of larger loads and can exercise more market power.

\section{Discussions}\label{sec_pre_6}
\textcolor{black}{
In this section, we discuss the limitations of the study and potential implications for policymakers.}

\subsection{Limitations of the study}
\textcolor{black}{The closed-form analysis of supply function equilibrium in a two-stage settlement market is theoretically complex and the system is often analyzed under certain simplifying market assumptions. While there are works that certainly consider a more relaxed set of assumptions, they either study single-stage market~\cite{Li_Na_linear_SFE,Joh2004Efficiency,yue_mpm, bansal2021market}, competitive market structure~\cite{mpm_mechanism_3, two_stage_competition_1,bansal2021market}, inelastic demand~\cite{johari_inelastic,Holmberg_inelastic}, homogeneous participants or symmetric participation~\cite{pyou_discovery,klemper_meyer}, etc. Given that the supply-function Nash equilibria are generally hard to characterize, even for a single-stage market, the literature often uses this approach as a zero-order analysis~\cite{pyou_discovery,klemper_meyer}. We further note that our findings are consistent with the theory in our numerical
experiments, where we relax, for e.g., the homogeneity constraints. Although capacity constraints, network constraints, etc., impact market power, our focus is on system-level market power, which occurs regardless of these constraints. Furthermore, an analysis
contemplating all of such settings, though interesting, will be too nuanced and digressed from our goal of counterfactual study.
}
%Most existing work, considers single-stage market~\cite{Li_Na_linear_SFE,Joh2004Efficiency, bansal2021market}, competitive market structure~\cite{mpm_mechanism_3, two_stage_competition_1}, symmetric participation~\cite{pyou_discovery}, etc. for discussion purposes. With the goal of a counterfactual study to understand the impact of system-level policy initiatives and immediate concerns of market power, we consider a simplified setting for the analysis.   The study concludes that DA-MPM policies are better suited for market mitigation. The future research direction will involve a more rigorous study of this policy, taking into account generator capacity and ramp constraints, network constraints, elastic demand, etc.}

\subsection{Implications for policymakers}
\textcolor{black}{%The evidence of (approximately $2\%$ hours) periods with non-competitive bids in spite of LMPM policies has led to the development of system-level initiatives. While arguably, in an ideal world, one would require the implementation of such changes (system-level market power mitigation) simultaneously in the day-ahead and the real-time market, the implementation of these changes is currently being considered independently. 
Our work highlights the importance of counterfactual analysis of two specific system-level policies in the CAISO area. Despite the CAISO's proposal of implementing a real-time MPM policy in the first phase, it should not be deployed by itself. We show that such a policy results in an inefficient competitive equilibrium, while the Nash equilibrium does not exist. We believe that if a strategy does not work well in a simple setting, then is unlikely to do well in a more complicated one. A day-ahead policy seems to have a reasonable impact on the market outcome that merits further analysis (with capacity constraints, network constraints, etc.) and consideration. Despite errors in cost estimations of generators, it results in an efficient competitive equilibrium, meaning the outcome aligns with the social planner problem. Moreover, the Nash equilibrium is more robust to market power and price manipulations. The aggregate profit (payment) of participants at Nash equilibrium is comparatively closer to the competitive one. Furthermore, the impact of error in cost estimation, heterogeneity in cost coefficient, and diversity in load size help policymakers with the tools to ensure fairness in the market. A positive bias in estimation leads to more profits for generators and mitigation of market power of loads. Larger loads may tend to split into smaller loads that merit further analysis to ensure market fairness.
}
%\textcolor{black}{The evidence of (approximately $2\%$ hours) periods with non-competitive bids in spite of LMPM policies has refocused the importance of understanding the participant behavior and investigation of competition in the market. It also led to the development of system-level initiatives with preemptive MPM policies, i.e., bid mitigation similar to LMPM, but at a system level for each stage before the market dispatch run. In particular, CAISO proposed an RT-MPM policy in the first phase arguing that real-time is more susceptible to market power and day-ahead is relatively competitive due to additional measures like virtual bidding that add competitive pressure on market clearing. However, our study highlights that even in a simple setting, a real-time MPM policy does not achieve efficient competitive equilibrium. Additionally, in the presence of strategic participants, a stable market outcome does not exist. Instead, a day-ahead MPM policy is more effective in mitigating market power to some extent.}

\section{Conclusions}\label{sec_6}

\textcolor{black}{We study competition between generators (bid linear supply function) and loads (bid quantity) in a two-stage settlement electricity market with a stage-wise MPM policy. In the proposed policy framework, CAISO substitutes generator bids with default bids in the stage with an MPM policy, i.e., day-ahead or real-time. To understand the participant behavior in the market, we start with a real-time MPM policy and analyze the sequential game, where generators only bid in the day-ahead market. The resulting competitive equilibrium, price-taker participants, is inefficient, while the Nash equilibrium, price-anticipating participants, does not exist, indicating an unstable market outcome. }

\textcolor{black}{Despite the estimation error, in the case of a day-ahead MPM policy, the competitive equilibrium aligns with the social planner problem. Further, the Nash equilibrium is robust to price manipulations compared to the standard market. Notably, our analysis shows that demand, despite being inelastic, could shift its allocation to manipulate market prices and win the competition. A more nuanced analysis of cost estimation error and heterogeneity in cost coefficients benefits generator over loads.} In the case of heterogeneous generators, expensive generators are less affected in the market. Also, the load size diversity highlights the role of a sufficiently smaller load in exercising market power at the expense of larger loads.

\bibliographystyle{IEEEtran}
\bibliography{TSG}

%\end{thebibliography}

\begin{IEEEbiography}[{\includegraphics[width=1in,height=1.25in,clip,keepaspectratio]{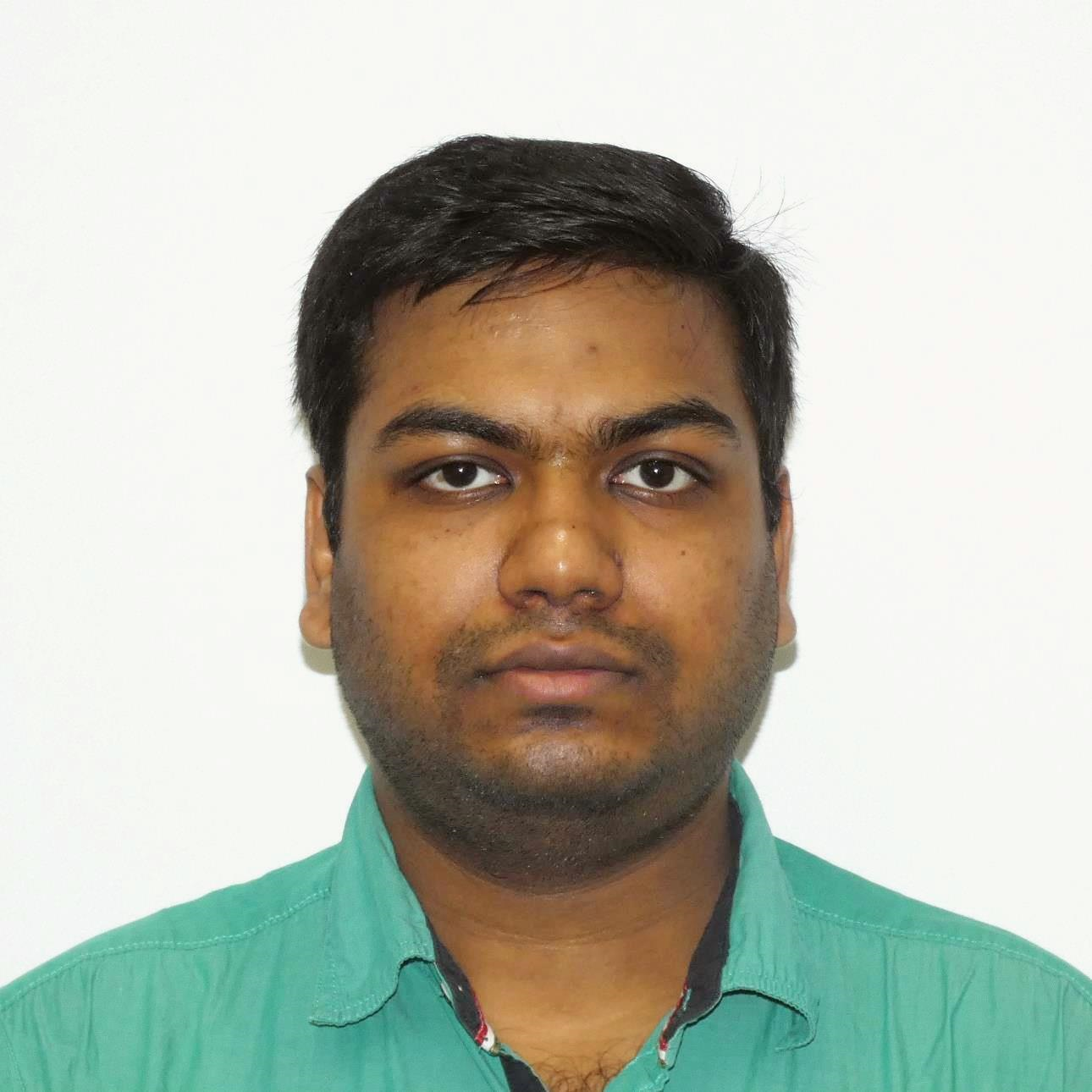}}]{Rajni Kant Bansal} received a B.Tech. in Mechanical Engineering from the Indian Institute of Technology Kanpur in 2016, where he received Academic Excellence Award. He was an
intern at John F. Welch Technology Centre, General Electric, in 2015. After graduating, he worked as an Analyst at Credit Suisse in India from 2016 to 2018. He then entered the
Ph.D. program in the Department of Mechanical Engineering in August 2018 and M.S.E program in Applied Mathematics and Statistics in 2022 at Johns Hopkins University. His
research uses traditional analytic models and tools from the operation research community, in particular graph theory, applied optimization, and mathematical economics, in order to
design mechanisms for efficient resource allocation in markets.
\end{IEEEbiography}

\begin{IEEEbiography}[{\includegraphics[width=1in,height=1.25in,clip,keepaspectratio]{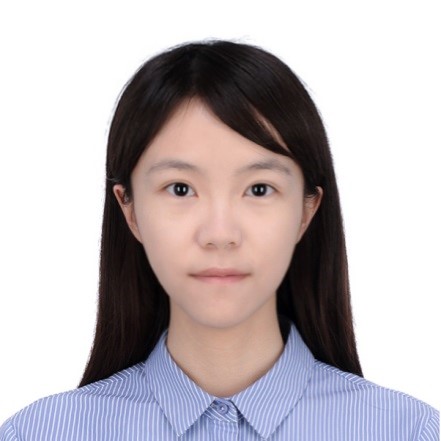}}]{Yue Chen} (Member, IEEE) received the B.E. degree in Electrical engineering from Tsinghua University, Beijing, China, in 2015, the B.S. degree in Economics from Peking University, Beijing, China, in 2017, and the Ph.D. degree in electrical engineering from Tsinghua University, in 2020. She is currently a Vice-Chancellor Assistant Professor with the Department of Mechanical and Automation Engineering, the Chinese University of Hong Kong, Hong Kong SAR. Her research interests include optimization, game theory, mathematical economics, and their applications in smart grids and integrated energy systems. She is an Associate Editor of IEEE Transactions on Smart Grid, IEEE Power Engineering Letters, and IET Renewable Power Generation.
\end{IEEEbiography}

\begin{IEEEbiography}[{\includegraphics[width=1in,height=1.25in,clip,keepaspectratio]{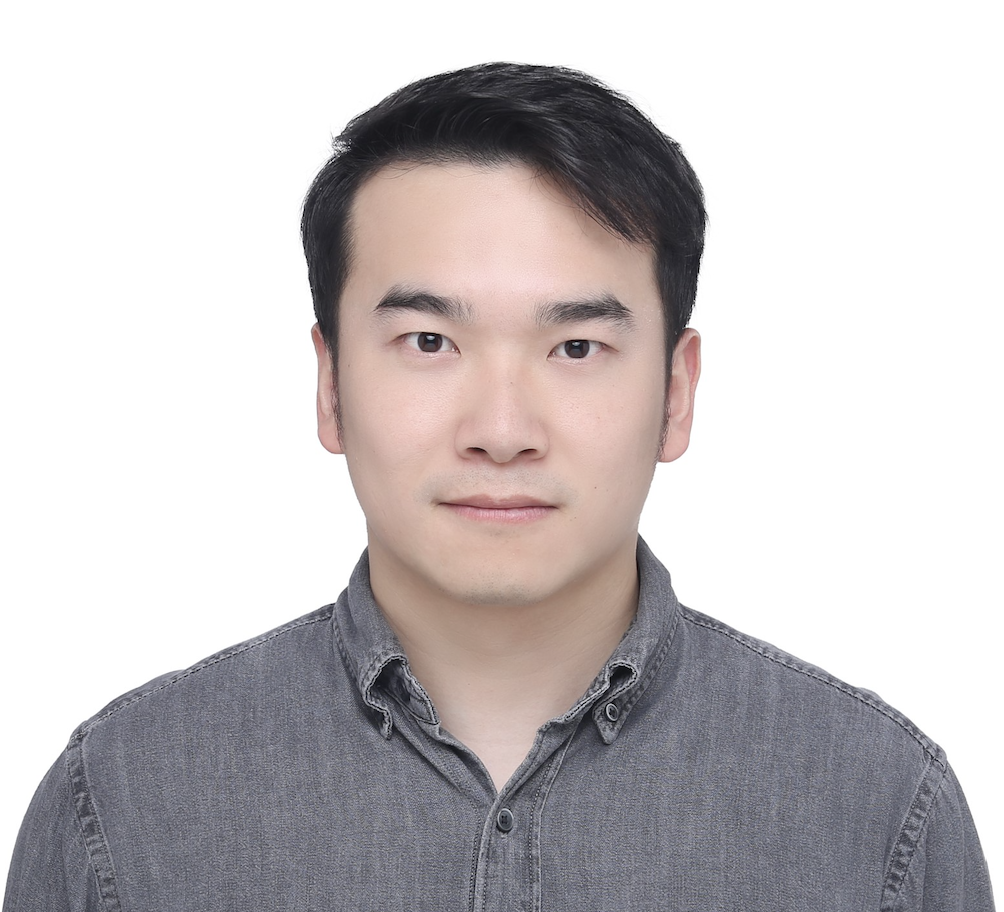}}]
{Pengcheng You} (S'14-M'18) is an Assistant Professor at the Department of Industrial Engineering and Management, Peking University. He also holds a joint appointment at the National Engineering Laboratory for Big Data Analysis and Applications, Peking University. Prior to joining PKU, he was a Postdoctoral Fellow at the ECE and ME Departments, Johns Hopkins University. He earned his Ph.D. and B.S. degrees both from Zhejiang University, China. During the graduate studies, he was a visiting student at Caltech and a research intern at PNNL.
His research interests include control, optimization, reinforcement learning and market mechanism, with main applications in power and energy systems.
\end{IEEEbiography}

\begin{IEEEbiography}[{\includegraphics[width=1in,height=1.25in,clip,keepaspectratio]{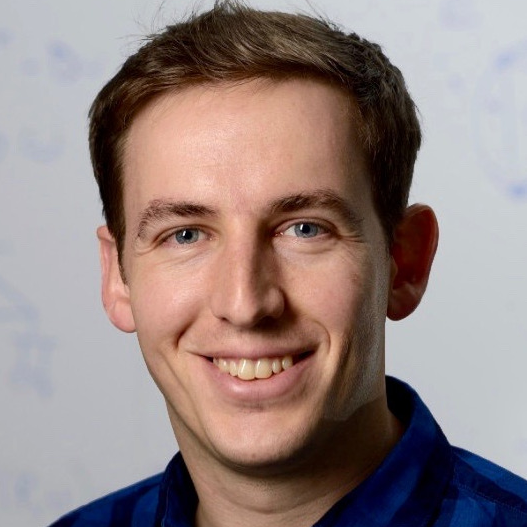}}]{Enrique Mallada} (S'09-M'13-SM') is an Associate Professor of Electrical and Computer Engineering at Johns Hopkins University. He was an Assistant Professor in the same department from 2016 to 2022. Prior to joining Hopkins, he was a Post-Doctoral Fellow in the Center for the Mathematics of Information at Caltech from 2014 to 2016. He received his Ingeniero en Telecomunicaciones degree from Universidad ORT, Uruguay, in 2005 and his Ph.D. degree in Electrical and Computer Engineering with a minor in Applied Mathematics from Cornell University in 2014.
Dr. Mallada was awarded
the NSF CAREER award in 2018,
the ECE Director's PhD Thesis Research Award for his dissertation in 2014,
the Center for the Mathematics of Information (CMI) Fellowship from Caltech in 2014,
and the Cornell University Jacobs Fellowship in 2011.
His research interests lie in the areas of control, dynamical systems, optimization, machine learning, with applications to engineering networks such as power systems and the Internet.
\end{IEEEbiography}

\ifthenelse{\boolean{arxiv}}{\color{black}
\appendices
\section{Proof of Theorem~\ref{comp_eqbm_rt_mpm_Thrm1}}\label{appendix_comp_eqbm_rt_mpm}

Under price-taking behavior, the individual problem for loads~\eqref{load_payment_final_rt_mpm} is a linear program with the closed-form solution given by:
\begin{align}\label{comp_eqbm_load_solution_rt_mpm}
    \!\!\!\!\left\{\begin{array}{l}
\!d_l^{d} \!=\!\infty, d_l^{r} \!=\! -\infty, d_l^{d}+d_l^{r} \!=\! d_l, \!\mbox{if } \lambda^{d} \!<\! \frac{d}{\sum\limits_{k\in\mathcal{G}}\!(c_k+\epsilon_k)^{-1}} \\
\!d_l^{d} \!=\! -\infty, d_l^{r} \!=\! \infty, d_l^{d}+d_l^{r} \!=\! d_l,  \mbox{if } \lambda^{d} \!>\! \frac{d}{\sum\limits_{k\in\mathcal{G}}\!(c_k+\epsilon_k)^{-1}}  \\
\!d_l^{d}+d_l^{r} = d_l, \quad  \mbox{ if } \lambda^{d} = \frac{d}{\sum\limits_{k\in\mathcal{G}}\!(c_k+\epsilon_k)^{-1}}
\end{array}\right.
\end{align}
where loads prefer the lower price in the market. The individual problem for generators~\eqref{generator_price_taking_profit_final_rt_mpm} requires:
\begin{align}\label{comp_eqbm_gen_solution_rt_mpm}
    \left\{\begin{array}{l}
\theta_j^d = -\infty, \!\mbox{ if } 0 \le \lambda^{d} < \frac{d}{\sum\limits_{k\in\mathcal{G}}\!(c_k+\epsilon_k)^{-1}} \\
\theta_j^d = \infty, \!\mbox{ if } \lambda^{d} < \frac{d}{\sum\limits_{k\in\mathcal{G}}\!(c_k+\epsilon_k)^{-1}}, \textrm{ and }  \lambda^{d} < 0\\
\theta_j^d = \infty,  \mbox{ if } \lambda^{d} > \frac{d}{\sum\limits_{k\in\mathcal{G}}\!(c_k+\epsilon_k)^{-1}}  \\
\theta_j^d \in \mathbb{R}_{\ge 0}, \quad  \mbox{ if } \lambda^{d} = \frac{d}{\sum\limits_{k\in\mathcal{G}}\!(c_k+\epsilon_k)^{-1}}
\end{array}\right.
\end{align}
where generators prefer higher prices in the market and seek to maximize profit. 
At the competitive equilibrium the day-ahead supply function~\eqref{gen_da_bid}, real-time true dispatch condition~\eqref{rt_true_dispatch}, real-time clearing prices~\eqref{rt_true_prc}, and the individual optimal solution \eqref{comp_eqbm_load_solution_rt_mpm},\eqref{comp_eqbm_gen_solution_rt_mpm} holds simultaneously and this is only possible if the market price is equal in the two-stages, i.e., 
\[
    \lambda^{d} = \lambda^{r} = \frac{d}{\sum\limits_{k\in\mathcal{G}}\!(c_k+\epsilon_k)^{-1}}, \textrm{s.t } d_l = d_l^d+d_l^r
\]
From real-time true dispatch conditions we have 
\[
    g_j^r+g_j^d = \frac{(c_j+\epsilon_j)^{-1}d}{\sum\limits_{k\in\mathcal{G}}\!(c_k+\epsilon_k)^{-1}}
\]
Thus a set of competitive equilibria exists.

\section{Proof of Theorem~\ref{strat_thrm_rt_mpm}}\label{appendix_strat_eqbm_rt_mpm}

From the day-ahead market clearing we have 
\begin{align}\label{day_ahead_mkt_clearing_rt_mpm}
    \sum_{j\in\mathcal{G}}\theta_j^d\lambda^d = d^d \implies \lambda^d = \frac{1}{\sum_{j\in\mathcal{G}}\theta_j^d}d^d
\end{align}
where we assume that $\sum_{j\in\mathcal{G}}\theta_j^d \neq 0$. Substituting~\eqref{day_ahead_mkt_clearing_rt_mpm} in generator individual profit optimization~\eqref{generator_strategic_profit_rt_mpm}, we get the individual problem of strategic generator $j$ as (we assume that $d^d\neq 0$ and leave the discussion of $d^d =0$ for later): 
\begin{align}
   \!\!\! \max_{\theta_j^d} & \left(\!\frac{d^d}{\sum\limits_{k \in \mathcal{G}}\theta_k^d} \!-\!\frac{d}{\sum\limits_{k \in \mathcal{G}}(c_k+\epsilon_k)^{-1}}\!\right)\!\frac{\theta_j^{d}d^d}{\sum\limits_{k \in \mathcal{G}}\theta_k^d}\! \nonumber\\
   &\quad \quad\quad \quad \quad\quad \quad  +\!  \frac{\frac{c_j}{2}+\epsilon_j}{(c_j+\epsilon_j)^2}\!\!\left(\!\frac{d}{\sum\limits_{k \in \mathcal{G}}(c_k+\epsilon_k)^{-1}}\!\!\right)^2
\end{align}
Though the individual problem is not necessarily concave in the domain, we can analyze the optimal bidding behavior from the first-order and second-order conditions. Writing the first-order condition, we have 

\begin{small}
\begin{align}
    & \frac{d \pi_j}{d \theta_j^d} =  \!\! \left[\!\frac{\theta_j^d}{\sum\limits_{k \in \mathcal{G}}\theta_k^d}\!\left(\frac{d}{\sum\limits_{k \in \mathcal{G}}(c_k+\epsilon_k)^{-1}}\!-\!\frac{2d^{d}}{\sum\limits_{k \in \mathcal{G}}\theta_k^d} \right)\!\! \right. \nonumber \\
    & \quad \quad\quad\quad\quad\quad\left.+ \!\!\left(\!\frac{d^{d}}{\sum\limits_{k\in\mathcal{G}}\theta_k^d}\!-\!\frac{d}{\sum\limits_{k \in \mathcal{G}}(c_k+\epsilon_k)^{-1}}\!\right)\!\!\right]\frac{d^{d}}{\sum\limits_{k \in \mathcal{G}}\theta_k^d} \label{strat_thrm_proof_rt_mpm_first_order}
\end{align}
\end{small}
Now summing over $j\in \mathcal{G}$ to attain the turning point of~\eqref{strat_thrm_proof_rt_mpm_first_order}, we have
\begin{small}
\begin{align}
    \implies & (G-2)(d^{d})^2 - (G-1)\frac{\sum_{j\in\mathcal{G}}\theta_j^d}{\sum_{k\in \mathcal{G}}(c_k+\epsilon_k)^{-1}}dd^{d} = 0
\end{align}
\end{small}
 where we assume that $|\mathcal{G}| \ge 2$. For the assumption $d^d \neq 0$, the potential turning point is given by
%\begin{align}\label{strat_thrm_rt_mpm_gen_opt_sol}
%    &\left\{\begin{array}{l} d^{d} = 0, \ \theta_j^d \in \mathbb{R}_{\ge 0}, \textrm{ or, } \\
%     d^{d} \neq 0, \theta_j^d = \frac{1}{G}\left(\sum_{j }c_j^{-1}\right)\frac{G-2}{G-1}\frac{d^{d}}{d}\end{array}\right.
%\end{align}

\begin{align}\label{strat_thrm_rt_mpm_gen_opt_sol}
    & \theta_j^d = \frac{1}{G}\left(\sum_{k \in \mathcal{G}}(c_k+\epsilon_k)^{-1}\right)\frac{G-2}{G-1}\frac{d^{d}}{d}
\end{align}

Similarly, substituting~\eqref{day_ahead_mkt_clearing_rt_mpm} in load individual payment optimization~\eqref{load_strategic_payment_rt_mpm}, we get the individual problem of load $l$ as - 
\begin{align}\label{strat_proof_load_optmz_rt_mpm}
    & \min_{d_l^{d}} \ \frac{d^{d}}{\sum_{j\in\mathcal{G}}\theta_j^d}d_l^{d} +\frac{d}{\sum_{j \in \mathcal{G}}(c_j+\epsilon_j)^{-1}}(d_l-d_l^{d})    
\end{align}
The unique optimal solution to the quadratic program~\eqref{strat_proof_load_optmz_rt_mpm} is given by 
\begin{align}\label{strat_thrm_rt_mpm_load_opt_sol}
    d_l^{d} = \frac{1}{L+1}\frac{\sum_{j\in \mathcal{G}}\theta_j^d}{\sum_{k \in \mathcal{G}}(c_k+\epsilon_k)^{-1}}d,  \ d_l^{r} = d_l - d_l^d
\end{align}
At equilibrium~\eqref{day_ahead_mkt_clearing_rt_mpm},\eqref{strat_thrm_rt_mpm_gen_opt_sol}, and \eqref{strat_thrm_rt_mpm_load_opt_sol} must hold simultaneously. This implies that
\[
    d^d = 0, \theta_j^d = 0 \implies \lambda^d = \lambda^r = \frac{1}{\sum_{j\in\mathcal{G}}{(c_j+\epsilon_j)^{-1}}}d
\]
where we use Rule~\ref{assumpt_same_prc} to define prices in the day-ahead market. However, this is in contradiction to our assumption and can be rejected. 

In the case of $d^d = 0$, 
\begin{itemize}
    \item If $\sum_{j\in\mathcal{G}} \theta_d \neq 0$, then solving ~\eqref{day_ahead_mkt_clearing_rt_mpm} and \eqref{strat_thrm_rt_mpm_load_opt_sol} simultaneously implies that $\sum_{j\in\mathcal{G}} \theta_d = 0$, which contradicts our assumption.
    \item If $\sum_{j\in\mathcal{G}} \theta_d = 0$, then we define prices using the Rule~\ref{assumpt_same_prc} in the day-ahead market. However, in this case, loads have the incentive to deviate from the equilibrium by allocating some demand in the day-ahead market since $\lambda^d = 0$, meaning loads make zero payment in the market, using Rule~\ref{assumpt_zero_prc}. 
\end{itemize}
Therefore the equilibrium does not exist. Similarly, in the case of only one generator, equilibrium does not exist. Though the generator bids arbitrary small values in the day ahead to earn increasing revenue, the load will also bid small quantities to decrease its payment. Since the generator operates truthfully in real-time, we attain the same equilibrium with all the demand allocated to the real-time market. Again, loads have the incentive to deviate and allocate demand in the day ahead where prices are zero. This completes the proof of Theorem~\ref{strat_thrm_rt_mpm}.

\section{Proof of Theorem~\ref{comp_eqbm_da_mpm_thrm}}\label{appendix_comp_eqbm_da_mpm}

Under price-taking behavior, the individual problem for loads~\eqref{load_payment_bids} is a linear program with the closed-form solution given by:
\begin{align}\label{comp_eqbm_load_solution_da_mpm}
    \!\!\!\!\left\{\begin{array}{l}
d_l^{d} = \infty, d_l^{r} = -\infty, d_l^{d}+d_l^{r} = d_l, \!\mbox{ if }  \lambda^d < \lambda^{r} \\
d_l^{d} = -\infty, d_l^{r} = \infty, d_l^{d}+d_l^{r} = d_l,  \mbox{ if } \lambda^d > \lambda^{r} \\
d_l^{d}+d_l^{r} = d_l, \quad  \mbox{ if } \lambda^d = \lambda^{r} 
\end{array}\right.
\end{align}
where loads prefer the lower price in the market. Similarly, solving concave individual problem of each generator~\eqref{generator_price_taking_profit_final_da_mpm} by taking the derivative, we have

%\begin{subequations}
\begin{small}
\begin{align}\label{comp_eqbm_gen_kkt_solution_da_mpm}
    %&{\lambda^r}^2 -c_j(\theta_j^r\lambda^{r}+\frac{c_j^{-1}d^d}{\sum_jc_j^{-1}})\lambda^{r} = 0 \nonumber \\\implies 
    &\!\!\!\!\lambda^{r}\!\left(\!\!(1\!-\!c_j\theta_j^r)\lambda^{r}\!-\!c_j\omega_jd^d\!\right)\!\! =\! 0 \!\!\implies \!\!\theta_j^r\lambda^r \!\!= \!c_j^{-1}\lambda^{r}\!\!-\!w_jd^d%\label{comp_eqbm_gen_kkt_solution_da_mpm}
\end{align}
\end{small}
%\end{subequations}

where we assume $\lambda^r \neq 0$. Summing~\eqref{comp_eqbm_gen_kkt_solution_da_mpm} over $j \in \mathcal{G}$ and using real-time market clearing~\eqref{rt_power_bal},
%\begin{align}\label{real_time_mkt_clearing_da_mpm}
%    \sum_{j\in\mathcal{G}}\theta_j^r\lambda^r = d^r
%\end{align}
we get 
\begin{align}
    %&{\lambda^r}^2 -c_j(\theta_j^r\lambda^{r}+\frac{c_j^{-1}d^d}{\sum_jc_j^{-1}})\lambda^{r} = 0 \nonumber \\\implies 
    &d^r= \sum_{j\in\mathcal{G}}c_j^{-1}\lambda^{r}-d^d \implies \lambda^r = \frac{d}{\sum_{j \in  \mathcal{G}}c_j^{-1}}\label{real_time_prc_da_mpm}
\end{align}
At equilibrium~\eqref{da_true_prc}, \eqref{comp_eqbm_load_solution_da_mpm}, \eqref{comp_eqbm_gen_kkt_solution_da_mpm}, and \eqref{real_time_prc_da_mpm} must hold simultaneously. This implies that  
\begin{align}
    & \lambda^d = \lambda^r = \frac{d}{\sum_{j\in\mathcal{G}}c_j^{-1}},\ d^d = \frac{\sum_{j\in\mathcal{G}}(c_j+\epsilon_j)^{-1}}{\sum_{j\in\mathcal{G}}c_j^{-1}}d\\
    & g_j^d = \omega_jd^d, \theta_j^r =\frac{\epsilon_j}{c_j(c_j+\epsilon_j)}, g_j^r = \frac{\epsilon_j}{c_j(c_j+\epsilon_j)}\lambda^r
\end{align}
Hence the equilibrium exists, and this completes the proof of Theorem~\ref{comp_eqbm_da_mpm_thrm}. 

\section{Proof of Theorem~\ref{strat_thrm_da_mpm}}\label{appendix_strat_eqbm_da_mpm}
Using the real-time clearing~\eqref{rt_power_bal}, we have 
\begin{align}\label{real_time_mkt_clearing_da_mpm_v2}
    \lambda^r = \frac{d^r}{\sum_{j\in\mathcal{G}}\theta_j^r}
\end{align}

where we assume that $\sum_{j\in\mathcal{G}}\theta_j^{r} \neq 0$. Substituting~\eqref{real_time_mkt_clearing_da_mpm_v2} in generator individual problem~\eqref{generator_strategic_profit_da_mpm}, the individual problem of price-anticipating generator $j$ is given by:

\begin{align}
    \!\!&\!\!\!\max_{\theta_j^r} (c_j+\epsilon_j)(\omega_jd^d)^2\! +\! \frac{\theta_j^r{d^r}^2}{(\sum\limits_{k\in\mathcal{G}}\theta_k^r)^2} \!-\!  \frac{c_j}{2}\!\!{\left(\!\!\omega_jd^d\! +\! \frac{\theta_j^rd^r}{\sum\limits_{k\in\mathcal{G}}\theta_k^r}\!\!\right)\!\!}^2 
\end{align}
We again use first-order and second-order conditions to analyze the optimal bidding behavior since the individual problem may not be concave in the domain. Writing the first-order condition, we have
\begin{align}\label{strat_thrm_proof_da_mpm_first_order_v1}
        \frac{d {\pi_j}}{d {\theta_j^r}} = 
        \frac{d^{r}}{(\sum\limits_{k\in\mathcal{G}}\theta_{k}^r)^3}\left[m_j^r - n_j^r \theta_j^r\right]
\end{align}
where 
\[
    m_j^r := d^{r}\!\!\sum_{k, k\neq j}\!\!\theta_k^r - \frac{c_j}{c_j+\epsilon_j}\frac{d^d}{\sum\limits_{k\in\mathcal{G}}(c_k+\epsilon_k)^{-1}}({\sum_{k, k\neq j}\!\!\theta_k^r})^2
\]
and 
\[
    n_j^r := d^{r} + \frac{c_j}{c_j+\epsilon_j}\frac{d^d}{\sum\limits_{k\in\mathcal{G}}(c_k+\epsilon_k)^{-1}}\sum_{k, k\neq j}\!\!\theta_k^r+c_jd^{r}\sum_{k, k\neq j}\!\!\theta_k^r
\]
Assuming generators are homogeneous and bid symmetrically, we can rewrite~\eqref{strat_thrm_proof_da_mpm_first_order_v1} as 
\begin{align}\label{strat_thrm_proof_da_mpm_first_order_v2}
        \frac{d \pi_j}{d \theta_j^r} = 
        \frac{d^{r}}{G^3{\theta^r}^2}\left[d^{r}(G-2)- cd(G-1)\theta^r\right]
\end{align}
then the turning point is given by 
\begin{align}\label{strat_proof_turning_pt_da_mpm}
    \theta^r_p = \frac{G-2}{G-1}\frac{d^r}{cd}
\end{align}
Writing the second-order condition and evaluating for homogeneous generators that bid symmetrically, i.e., the turning point~\eqref{strat_proof_turning_pt_da_mpm}, we have
\begin{align}\label{strat_thrm_proof_da_mpm_second_order_v1}
        &\frac{d^2 {\pi_j}}{d {\theta_j^r}^2}\bigg|_{\theta_j^r = \theta^r_p(d^r)} = 
        \frac{d^{r}}{(\sum_{j}\theta_{j}^r)^3}\left[\tilde{m}_j^r + \tilde{n}_j^r \theta_j^r\right]\bigg|_{\theta_j^r = \theta^r_p(d^r)} \\
        & \quad \quad = -\frac{c^3(G-1)^4d^2}{G^4(G-2)^3}\frac{d}{d^r}\left(2+ (G-2)\frac{d^r}{d}\right)
\end{align}
where 
\[
    \tilde{m}_j^r := -\frac{4d^{r}\sum_{k, k\neq j}\theta_k^r}{\sum_k\theta_k^r} + 2c_j\omega_jd^d{\sum_{k, k\neq j}\theta_k^r} - \frac{c_jd^{r}(\sum_{k, k\neq j}\theta_k^r)^2}{\sum_k\theta_k^r}
\]
and 
\[
    \tilde{n}_j^r := \frac{2d^{r}}{\sum_k\theta_k^r} + \frac{2c_jd^{r}\sum_{k, k\neq j}\theta_k^r}{\sum_k\theta_k^r} 
\]
Now, loads acting as leaders anticipate the clearing prices and optimal bids of generators in the real-time subgame equilibrium, such that 
\begin{align}\label{strat_thrm_proof_da_mpm_real_time_prc}
    \lambda^r = \frac{G-1}{G-2}\frac{cd}{G}
\end{align}
where we substitute~\eqref{strat_proof_turning_pt_da_mpm} in \eqref{real_time_mkt_clearing_da_mpm_v2}. Substituting~\eqref{strat_thrm_proof_da_mpm_real_time_prc} in load individual problem~\eqref{load_strategic_payment_da_mpm}, we have 
\begin{align}\label{strat_proof_load_optmz_da_mpm}
    & \min_{d_l^{d}} \frac{d^d}{\sum_{j\in\mathcal{G}}(c_j+\epsilon_j)^{-1}}d_l^{d} + \frac{G-1}{G-2}\frac{cd}{G}(d_l-d_l^{d})
\end{align}
The unique optimal solution to the quadratic program~\eqref{strat_proof_load_optmz_da_mpm} is given by
\begin{align}
    d_l^{d} = \frac{c}{c+\epsilon}\frac{1}{L+1}\frac{G-1}{G-2}d,  \ d_l^{r} =\left(d_l - \frac{c}{c+\epsilon}\frac{1}{L+1}\frac{G-1}{G-2}d\right)
\end{align}
where we assume that generators are homogeneous and estimation error is the same, i.e. i.e., $c_j:=c, \epsilon_j: = \epsilon \ \forall j \in \mathcal{G}$ 
Assuming $\frac{1}{L} > \frac{c-\epsilon(G-2)}{(c+\epsilon)(G-2)}$,
\[
    d^r > 0 \implies \frac{d^2 {\pi_j}}{d {\theta_j^r}^2}\bigg|_{\theta_j^r = \theta^r_p} < 0
\]
Thus the obtained equilibrium maximizes generators' profit and minimizes loads' payment while the supply-demand balance is satisfied. However, if $\frac{1}{L} < \frac{c-\epsilon(G-2)}{(c+\epsilon)(G-2)}$, then 
\[
    d^r < 0 \implies \theta^r_p < 0 \implies \frac{d^2 {\pi_j}}{d {\theta_j^r}^2}\bigg|_{\theta_j^r = \theta^r_p} > 0
\]
The obtained equilibrium minimizes generators' profit, and generators' have the incentive to deviate from this equilibrium. Therefore, symmetric equilibrium does not exist in this case. Moreover, in the case of $|\mathcal{G}|<3$, generators have the incentive to bid arbitrarily small values and earn arbitrarily large profits in the market.

In the case of $\frac{1}{L} = \frac{c-\epsilon(G-2)}{(c+\epsilon)(G-2)}$, at equilibrium $d^r =0$ which contradicts our initial assumption. We analyze the case $d^r=0$ separately, 
\begin{enumerate}
    \item If $\sum_{j\in \mathcal{G}}\theta_j^r \neq 0 \implies \lambda^r = 0$ and $\lambda^d = \frac{c+\epsilon}{G}d^d = \frac{c+\epsilon}{G}d$. WLOG, we can assume that $d_l^r=0, \ \forall l \in \mathcal{L}$, otherwise load $l$ with non-zero demand has the incentive to deviate and participate in the real-time market to minimize its payment.  The payment of individual load $l$ is then given by 
    \[
       \lambda^d d_l^d+\lambda^r d_l^r = \frac{c+\epsilon}{G}dd_l^d
    \]
    However, if load $l$ unilaterally decides to deviate by allocating demand in real-time, i.e., $d_l^r = \gamma$ then the payment is given by 
    \[
       \lambda^d d_l^d+\lambda^r d_l^r = \frac{c+\epsilon}{G}(d_l-\gamma)d + \frac{\gamma^2}{\sum_{j\in \mathcal{G}}\theta_j}
    \]
    which is smaller for small enough $\gamma$. Therefore the equilibrium does not exist. 
    \item If $\sum_{j\in \mathcal{G}}\theta_j^r = 0$, using Rule~\ref{assumpt_same_prc} we have $\lambda^r = \lambda^d$ and $\lambda^d = \frac{c+\epsilon}{G}d^d = \frac{c+\epsilon}{G}d$. However, if load $l$ unilaterally decides to deviate by allocating demand in real-time i.e., $d_l^r = \gamma$ then using Rule~\ref{assumpt_zero_prc} $\lambda^r = 0$. Therefore load has the incentive to deviate and allocate demand in the real-time market with zero clearing price. Hence equilibrium does not exist.
    
\end{enumerate}
This completes the proof of the Theorem~\ref{strat_thrm_da_mpm}. 
}{}

\end{document}